\documentclass[a4paper,pdftex,12pt]{article}
\usepackage{tikz}
\usepackage{caption}
\usepackage{subcaption}
\usepackage{times}
\usepackage{amsmath}
\usepackage{amsfonts}
\usepackage{amssymb}
\usepackage[linkcolor=black,urlcolor=black,citecolor=black]{hyperref}
\usepackage[all]{hypcap}
\hypersetup{linktocpage,colorlinks=true,pdfborder={0 0 0}}
\usepackage{geometry}
\title{Aperiodic Sets of Prototiles Extracted From the Penrose Rhomb Tiling}

\author{Mike Winkler\footnote{http://www.mikematics.de/}}

\date{\small Fakult\"at f\"ur Mathematik, Ruhr-Universit\"at Bochum, Germany, mike.winkler@ruhr-uni-bochum.de}

\usepackage{titling}
\begin{document}

  \maketitle\thispagestyle{empty}
  
  \begin{abstract}
    We present aperiodic sets of prototiles whose shapes are based on the well-known Penrose rhomb tiling. Some decorated prototiles lead to an exact Penrose rhomb tiling without any matching rules. We also give an approximate solution to an aperiodic monotile that tessellates the plane (including five types of gaps) only in a nonperiodic way.
  \end{abstract}

  \section{Introduction} 
  
  The Penrose P2 (kite and dart) and P3 (rhomb) tilings are certainly the most popular examples of nonperiodic tilings. Both tilings are strongly related and generate the same mld-class\footnote{Two tilings are called mld (mutually locally derivable), if one is obtained from the other in a unique way by local rules, and vice versa \cite{Bielefeld}.}. For more details we refer the reader to \cite{Bielefeld}\cite{Gruenbaum}\cite{Penrose}. In this article we will take a closer look at cutouts from the P3 tiling, which themselves form aperiodic sets of two prototiles (decorated or undecorated) and which have not yet been published or mentioned. The smallest structures of edge to edge connected rhombs we have found are given in Figure 1, whereby $T_1$ is not a connected tile, and $T_3$ represents a connected pair of $T_1$ and $T_2$. The aperiodic sets of $T_1,T_2$ or $T_2,T_3$ inevitably lead to an exact P3 tiling without matching rules. A slight modification of $T_1$ and $T_2$ will give us an aperiodic set of two undecorated tiles in the shape of a snake and a dog (Fig. 4).
  
  \begin{figure}[!ht]
  	\centering
  	\begin{tikzpicture}
  		[y=-1.0pt, x=1.0pt, scale=0.8]
  		\draw[line width=0.5pt]
  		(11.7258,601.3563) -- (15.2210,590.5971)
  		(15.2210,590.5971) -- (18.7170,601.3563)
  		(18.7170,601.3563) -- (30.0306,601.3563)
  		(15.2210,590.5971) -- (26.5346,590.5971)
  		(15.2210,590.5971) -- (18.7170,579.8375)
  		(30.0306,636.1747) -- (26.5346,625.4155)
  		(39.1826,608.0059) -- (30.0306,601.3563)
  		(24.3738,653.5835) -- (27.8698,664.3435)
  		(24.3738,653.5835) -- (35.6874,653.5835)
  		(24.3738,653.5835) -- (15.2210,646.9339)
  		(27.8698,664.3435) -- (18.7170,657.6931)
  		(18.7170,657.6931) -- (15.2210,646.9339)
  		(9.5650,664.3435) -- (18.7170,657.6931)
  		(35.6874,597.2467) -- (26.5346,590.5971)
  		(39.1826,608.0059) -- (35.6874,597.2467)
  		(35.6874,597.2467) -- (44.8394,590.5971)
  		(44.8394,590.5971) -- (35.6874,583.9471)
  		(35.6874,583.9471) -- (26.5346,590.5971)
  		(11.7258,601.3563) -- (15.2210,612.1155)
  		(11.7258,601.3563) -- (2.5730,608.0059)
  		(2.5730,608.0059) -- (6.0690,618.7651)
  		(24.3738,618.7651) -- (15.2210,625.4155)
  		(15.2210,612.1155) -- (24.3738,618.7651)
  		(15.2210,612.1155) -- (6.0690,618.7651)
  		(15.2210,625.4155) -- (6.0690,618.7651)
  		(15.2210,625.4155) -- (26.5346,625.4155)
  		(15.2210,625.4155) -- (18.7170,636.1747)
  		(35.6874,653.5835) -- (26.5346,646.9339)
  		(30.0306,636.1747) -- (18.7170,636.1747)
  		(18.7170,636.1747) -- (15.2210,646.9339)
  		(26.5346,646.9339) -- (15.2210,646.9339)
  		(30.0306,636.1747) -- (26.5346,646.9339)
  		(30.0306,601.3563) -- (26.5346,590.5971)
  		(15.2210,612.1155) -- (18.7170,601.3563)
  		(18.7170,579.8375) -- (30.0306,579.8375)
  		(30.0306,579.8375) -- (26.5346,590.5971)
  		(27.8698,664.3435) -- (18.7170,670.9931)
  		(9.5650,664.3435) -- (18.7170,670.9931)
  		(15.2210,625.4155) -- (11.7258,636.1747)
  		(11.7258,636.1747) -- (15.2210,646.9339)
  		(116.9541,646.9174) -- (113.4581,636.1574)
  		(113.4581,636.1574) -- (109.9621,625.3982)
  		(104.3053,664.3262) -- (92.9925,664.3262)
  		(102.1445,657.6766) -- (92.9925,664.3262)
  		(113.4581,636.1574) -- (109.9621,646.9174)
  		(109.9621,646.9174) -- (98.6485,646.9174)
  		(98.6485,646.9174) -- (102.1445,657.6766)
  		(113.4581,636.1574) -- (102.1445,636.1574)
  		(92.9925,664.3262) -- (89.4965,653.5670)
  		(89.4965,653.5670) -- (98.6485,646.9174)
  		(89.4965,653.5670) -- (78.1829,653.5670)
  		(78.1829,653.5670) -- (81.6789,664.3262)
  		(92.9925,642.8078) -- (102.1445,636.1574)
  		(89.4965,653.5670) -- (92.9925,642.8078)
  		(92.9925,642.8078) -- (83.8397,636.1574)
  		(83.8397,636.1574) -- (80.3437,646.9174)
  		(80.3437,646.9174) -- (89.4965,653.5670)
  		(83.8397,636.1574) -- (92.9925,629.5078)
  		(92.9925,629.5078) -- (102.1445,636.1574)
  		(92.9925,629.5078) -- (89.4965,618.7486)
  		(98.6485,625.3982) -- (89.4965,618.7486)
  		(116.9541,646.9174) -- (113.4581,657.6766)
  		(113.4581,657.6766) -- (104.3053,664.3262)
  		(104.3053,664.3262) -- (113.4581,670.9758)
  		(113.4581,657.6766) -- (122.6101,664.3262)
  		(113.4581,670.9758) -- (122.6101,664.3262)
  		(109.9621,625.3982) -- (98.6485,625.3982)
  		(98.6485,625.3982) -- (102.1445,636.1574)
  		(98.6485,646.9174) -- (102.1445,636.1574)
  		(113.4581,657.6766) -- (102.1445,657.6766)
  		(113.4581,657.6766) -- (109.9621,646.9174)
  		(92.9925,664.3262) -- (81.6789,664.3262)
  		(109.9621,625.3982) -- (100.8093,618.7486)
  		(100.8093,618.7486) -- (89.4965,618.7486)
  		(100.8093,618.7486) -- (104.3053,607.9894)
  		(104.3053,607.9894) -- (113.4581,601.3398)
  		(104.3053,607.9894) -- (92.9925,607.9894)
  		(92.9925,607.9894) -- (89.4965,618.7486)
  		(92.9925,607.9894) -- (102.1445,601.3398)
  		(102.1445,601.3398) -- (113.4581,601.3398)
  		(102.1445,601.3398) -- (98.6485,590.5798)
  		(80.3437,625.3982) -- (89.4965,618.7486)
  		(78.1829,618.7486) -- (89.4965,618.7486)
  		(69.0309,612.0990) -- (78.1829,618.7486)
  		(69.0309,612.0990) -- (80.3437,612.0990)
  		(80.3437,612.0990) -- (89.4965,618.7486)
  		(80.3437,612.0990) -- (83.8397,601.3398)
  		(83.8397,601.3398) -- (92.9925,607.9894)
  		(92.9925,607.9894) -- (89.4965,597.2302)
  		(80.3437,625.3982) -- (69.0309,625.3982)
  		(69.0309,625.3982) -- (78.1829,618.7486)
  		(69.0309,612.0990) -- (72.5261,601.3398)
  		(72.5261,601.3398) -- (83.8397,601.3398)
  		(89.4965,597.2302) -- (98.6485,590.5798)
  		(83.8397,636.1574) -- (80.3437,625.3982)
  		(151.1436,601.2809) -- (154.6396,590.5213)
  		(154.6396,590.5213) -- (158.1356,601.2809)
  		(154.6396,590.5213) -- (165.9524,590.5213)
  		(154.6396,590.5213) -- (158.1356,579.7619)
  		(163.7916,653.5081) -- (154.6396,646.8585)
  		(167.2876,664.2673) -- (158.1356,657.6177)
  		(158.1356,657.6177) -- (154.6396,646.8585)
  		(148.9828,664.2673) -- (158.1356,657.6177)
  		(175.1052,597.1713) -- (165.9524,590.5213)
  		(184.2572,590.5213) -- (175.1052,583.8716)
  		(175.1052,583.8716) -- (165.9524,590.5213)
  		(151.1436,601.2809) -- (154.6396,612.0401)
  		(151.1436,601.2809) -- (141.9908,607.9305)
  		(141.9908,607.9305) -- (145.4868,618.6897)
  		(154.6396,612.0401) -- (145.4868,618.6897)
  		(154.6396,625.3393) -- (145.4868,618.6897)
  		(154.6396,625.3393) -- (158.1356,636.0993)
  		(169.4484,636.0993) -- (158.1356,636.0993)
  		(158.1356,636.0993) -- (154.6396,646.8585)
  		(165.9524,646.8585) -- (154.6396,646.8585)
  		(169.4484,601.2809) -- (165.9524,590.5213)
  		(158.1356,579.7619) -- (169.4484,579.7619)
  		(169.4484,579.7619) -- (165.9524,590.5213)
  		(167.2876,664.2673) -- (158.1356,670.9169)
  		(148.9828,664.2673) -- (158.1356,670.9169)
  		(154.6396,625.3393) -- (151.1436,636.0993)
  		(151.1436,636.0993) -- (154.6396,646.8585)
  		(202.5628,646.8585) -- (199.0668,636.0993)
  		(199.0668,636.0993) -- (195.5708,625.3393)
  		(189.9140,664.2673) -- (178.6012,664.2673)
  		(187.7532,657.6177) -- (178.6012,664.2673)
  		(199.0668,636.0993) -- (195.5708,646.8585)
  		(195.5708,646.8585) -- (184.2572,646.8585)
  		(184.2572,646.8585) -- (187.7532,657.6177)
  		(199.0668,636.0993) -- (187.7532,636.0993)
  		(178.6012,664.2673) -- (175.1052,653.5081)
  		(175.1052,653.5081) -- (184.2572,646.8585)
  		(175.1052,653.5081) -- (163.7916,653.5081)
  		(163.7916,653.5081) -- (167.2876,664.2673)
  		(178.6012,642.7489) -- (187.7532,636.0993)
  		(175.1052,653.5081) -- (178.6012,642.7489)
  		(178.6012,642.7489) -- (169.4484,636.0993)
  		(169.4484,636.0993) -- (165.9524,646.8585)
  		(165.9524,646.8585) -- (175.1052,653.5081)
  		(169.4484,636.0993) -- (178.6012,629.4497)
  		(178.6012,629.4497) -- (187.7532,636.0993)
  		(178.6012,629.4497) -- (175.1052,618.6897)
  		(184.2572,625.3393) -- (175.1052,618.6897)
  		(202.5628,646.8585) -- (199.0668,657.6177)
  		(199.0668,657.6177) -- (189.9140,664.2673)
  		(189.9140,664.2673) -- (199.0668,670.9169)
  		(199.0668,657.6177) -- (208.2188,664.2673)
  		(199.0668,670.9169) -- (208.2188,664.2673)
  		(195.5708,625.3393) -- (184.2572,625.3393)
  		(184.2572,625.3393) -- (187.7532,636.0993)
  		(184.2572,646.8585) -- (187.7532,636.0993)
  		(199.0668,657.6177) -- (187.7532,657.6177)
  		(199.0668,657.6177) -- (195.5708,646.8585)
  		(178.6012,664.2673) -- (167.2876,664.2673)
  		(195.5708,625.3393) -- (186.4180,618.6897)
  		(186.4180,618.6897) -- (175.1052,618.6897)
  		(186.4180,618.6897) -- (189.9140,607.9305)
  		(189.9140,607.9305) -- (199.0668,601.2809)
  		(189.9140,607.9305) -- (178.6012,607.9305)
  		(178.6012,607.9305) -- (175.1052,618.6897)
  		(178.6012,607.9305) -- (187.7532,601.2809)
  		(187.7532,601.2809) -- (199.0668,601.2809)
  		(187.7532,601.2809) -- (184.2572,590.5213)
  		(165.9524,625.3393) -- (175.1052,618.6897)
  		(163.7916,618.6897) -- (175.1052,618.6897)
  		(154.6396,612.0401) -- (163.7916,618.6897)
  		(154.6396,612.0401) -- (165.9524,612.0401)
  		(165.9524,612.0401) -- (175.1052,618.6897)
  		(165.9524,612.0401) -- (169.4484,601.2809)
  		(169.4484,601.2809) -- (178.6012,607.9305)
  		(178.6012,607.9305) -- (175.1052,597.1713)
  		(165.9524,625.3393) -- (154.6396,625.3393)
  		(154.6396,625.3393) -- (163.7916,618.6897)
  		(154.6396,612.0401) -- (158.1356,601.2809)
  		(158.1356,601.2809) -- (169.4484,601.2809)
  		(175.1052,597.1713) -- (184.2572,590.5213)
  		(169.4484,636.0993) -- (165.9524,625.3393);
  		
  		\small
  		\coordinate[label=right:$T_1$] (a) at (15,685);
  		\coordinate[label=right:$T_2$] (a) at (90,685);
  		\coordinate[label=right:$T_3$] (a) at (170,685);
  		\normalsize
  	\end{tikzpicture}	
  	\caption{Cutouts of the P3 tiling.}
  \end{figure}
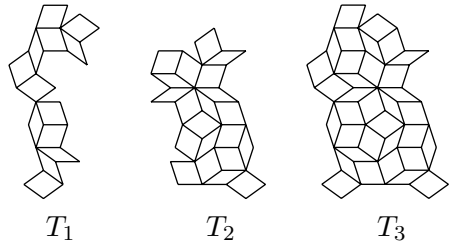
  
  Tilings like P2 and P3 have a scaling self-similarity as a fractal. This is due to the substitution rules that allow each tile to be decomposed into smaller tiles of the same shape as those used in the tiling. Thus allow larger tiles to be composed of smaller ones. Penrose tilings are strongly related with Fibonacci numbers\footnote{https://oeis.org/A000045}. For instance, the ratio between the number of different prototiles which can tessellate a larger tile is always given by such numbers \cite{Bielefeld}\cite{Gruenbaum}. This Fibonacci ratio also applies to the prototiles in Figure 1. $T_1$ consists of five thin and eight thick rhombs, $T_2$ of eight thin and thirteen thick rhombs, and $T_3$ of thirteen thin and 21 thick rhombs.
  
  Figures 2 and 3 show the arrangements for $T_1,T_2$ respectively $T_2,T_3$ based on the superordinate kite and dart shape. The use of the substitution rules for the P2 tiling then lead to a nonperiodic tiling (Fig. 10). Please note that the right sided arrangements of Figures 2 and 3 are completely part of the left sided ones.
  
  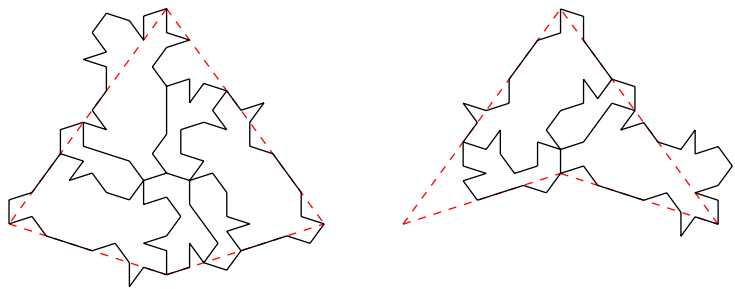
\begin{figure}[!ht]
  	\centering
  	\begin{tikzpicture}
  		[y=-1.0pt, x=1.0pt, scale=0.8]
  		\draw[red, dashed, line width=0.5pt]
  		(79.4446,539.0250) -- (5.6982,640.5280)
  		(5.6982,640.5280) -- (79.4446,664.4896)
  		(79.4446,664.4896) -- (153.1910,640.5281)
  		(79.4446,539.0250) -- (153.1910,640.5281)
  		(337.5449,640.6809) -- (263.7993,539.1780)
  		(263.7993,616.7194) -- (337.5449,640.6809)
  		(190.0529,640.6810) -- (263.7993,539.1780)
  		(190.0529,640.6810) -- (263.7985,616.7194);
  		
  		\draw[line width=0.5pt]
  		(90.2038,631.3755) -- (96.8534,640.5281)
  		(96.8534,640.5281) -- (103.5030,649.6806)
  		(40.5166,603.9180) -- (51.2758,607.4139)
  		(40.5166,592.6049) -- (40.5166,603.9180)
  		(79.4446,664.4896) -- (90.2038,660.9937)
  		(90.2038,660.9937) -- (90.2038,649.6806)
  		(90.2038,631.3755) -- (90.2038,620.0624)
  		(90.2038,620.0624) -- (79.4446,616.5664)
  		(103.5030,649.6806) -- (96.8534,658.8330)
  		(96.8534,658.8330) -- (90.2038,649.6806)
  		(29.7574,607.4139) -- (29.7574,596.1008)
  		(40.5166,592.6049) -- (29.7574,596.1008)
  		(23.1070,646.1846) -- (33.8670,649.6806)
  		(33.8670,649.6806) -- (44.6262,653.1765)
  		(16.4574,625.7190) -- (23.1070,616.5664)
  		(33.8670,620.0624) -- (40.5166,610.9099)
  		(40.5166,610.9099) -- (29.7574,607.4139)
  		(51.2758,625.7190) -- (44.6262,616.5664)
  		(44.6262,616.5664) -- (33.8670,620.0624)
  		(23.1070,646.1846) -- (16.4574,637.0321)
  		(16.4574,625.7190) -- (5.6982,629.2149)
  		(16.4574,637.0321) -- (5.6982,640.5280)
  		(5.6982,629.2149) -- (5.6982,640.5280)
  		(23.1070,616.5664) -- (29.7574,607.4139)
  		(44.6262,653.1765) -- (55.3854,649.6806)
  		(55.3854,649.6806) -- (62.0350,658.8330)
  		(62.0350,658.8330) -- (62.0350,670.1462)
  		(68.6854,660.9937) -- (62.0350,670.1462)
  		(68.6854,660.9937) -- (79.4446,664.4896)
  		(79.4446,627.8795) -- (68.6854,631.3755)
  		(68.6854,620.0624) -- (68.6854,631.3755)
  		(79.4446,627.8795) -- (86.0942,637.0321)
  		(86.0942,637.0321) -- (79.4446,646.1846)
  		(79.4446,598.2614) -- (79.4446,586.9482)
  		(51.2758,596.1008) -- (40.5166,592.6049)
  		(51.2758,577.7958) -- (44.6262,586.9482)
  		(44.6262,586.9482) -- (51.2758,596.1008)
  		(79.4446,598.2614) -- (72.7950,607.4139)
  		(72.7950,607.4139) -- (79.4446,616.5664)
  		(79.4446,575.6351) -- (72.7950,566.4826)
  		(72.7950,566.4826) -- (79.4446,557.3301)
  		(79.4446,557.3301) -- (90.2038,553.8342)
  		(79.4446,550.3382) -- (90.2038,553.8342)
  		(79.4446,550.3382) -- (79.4446,539.0250)
  		(44.6262,550.3382) -- (51.2758,559.4907)
  		(62.0350,544.6817) -- (68.6854,553.8342)
  		(68.6854,553.8342) -- (68.6854,542.5210)
  		(51.2758,566.4826) -- (40.5166,562.9867)
  		(40.5166,562.9867) -- (51.2758,559.4907)
  		(44.6262,550.3382) -- (51.2758,541.1857)
  		(51.2758,541.1857) -- (62.0350,544.6817)
  		(68.6854,542.5210) -- (79.4446,539.0250)
  		(51.2758,577.7958) -- (51.2758,566.4826)
  		(107.6134,577.7958) -- (96.8534,574.2998)
  		(96.8534,574.2998) -- (90.2038,583.4523)
  		(79.4446,575.6351) -- (90.2038,572.1392)
  		(90.2038,572.1392) -- (90.2038,583.4523)
  		(114.2630,653.1765) -- (107.6134,662.3290)
  		(96.8534,658.8330) -- (107.6134,662.3290)
  		(142.4310,625.7190) -- (135.7814,616.5664)
  		(135.7814,616.5664) -- (129.1318,607.4139)
  		(135.7814,646.1846) -- (125.0222,649.6806)
  		(118.3726,640.5281) -- (107.6134,644.0240)
  		(107.6134,644.0240) -- (114.2630,653.1765)
  		(107.6134,625.7190) -- (107.6134,637.0321)
  		(107.6134,637.0321) -- (118.3726,640.5281)
  		(142.4310,625.7190) -- (142.4310,637.0321)
  		(135.7814,646.1846) -- (146.5406,649.6806)
  		(142.4310,637.0321) -- (153.1910,640.5281)
  		(146.5406,649.6806) -- (153.1910,640.5281)
  		(125.0222,649.6806) -- (114.2630,653.1765)
  		(129.1318,607.4139) -- (118.3726,603.9180)
  		(118.3726,603.9180) -- (118.3726,592.6049)
  		(118.3726,592.6049) -- (125.0222,583.4523)
  		(114.2630,586.9482) -- (125.0222,583.4523)
  		(114.2630,586.9482) -- (107.6134,577.7958)
  		(86.0942,607.4139) -- (96.8534,610.9099)
  		(96.8534,592.6049) -- (107.6134,596.1008)
  		(107.6134,596.1008) -- (100.9630,586.9482)
  		(100.9630,616.5664) -- (90.2038,620.0624)
  		(90.2038,620.0624) -- (96.8534,610.9099)
  		(86.0942,607.4139) -- (86.0942,596.1008)
  		(86.0942,596.1008) -- (96.8534,592.6049)
  		(100.9630,586.9482) -- (107.6134,577.7958)
  		(107.6134,625.7190) -- (100.9630,616.5664)
  		(79.4446,616.5664) -- (68.6854,620.0624)
  		(79.4446,664.4896) -- (79.4446,653.1765)
  		(79.4446,653.1765) -- (68.6854,649.6806)
  		(68.6854,649.6806) -- (79.4446,646.1846)
  		(68.6854,620.0624) -- (62.0350,629.2149)
  		(62.0350,629.2149) -- (51.2758,625.7190)
  		(79.4446,575.6351) -- (79.4446,586.9482)
  		(68.6854,620.0624) -- (62.0350,610.9098)
  		(51.2758,607.4139) -- (62.0350,610.9098)
  		(246.3897,622.3760) -- (235.6305,625.8719)
  		(235.6305,625.8719) -- (224.8713,629.3678)
  		(287.8577,583.6052) -- (281.2081,592.7577)
  		(298.6169,587.1012) -- (287.8577,583.6052)
  		(218.2209,601.9103) -- (218.2209,613.2234)
  		(218.2209,613.2234) -- (228.9809,616.7194)
  		(246.3897,622.3760) -- (257.1489,625.8719)
  		(257.1489,625.8719) -- (263.7993,616.7194)
  		(224.8713,629.3678) -- (218.2209,620.2153)
  		(218.2209,620.2153) -- (228.9809,616.7194)
  		(287.8577,572.2921) -- (298.6169,575.7881)
  		(298.6169,587.1012) -- (298.6169,575.7881)
  		(253.0393,553.9871) -- (246.3897,563.1396)
  		(246.3897,563.1396) -- (239.7401,572.2921)
  		(274.5585,553.9871) -- (281.2081,563.1396)
  		(274.5585,572.2921) -- (281.2081,581.4446)
  		(281.2081,581.4446) -- (287.8577,572.2921)
  		(263.7993,587.1012) -- (274.5585,583.6052)
  		(274.5585,583.6052) -- (274.5585,572.2921)
  		(253.0393,553.9871) -- (263.7993,550.4912)
  		(274.5585,553.9871) -- (274.5585,542.6740)
  		(263.7993,550.4912) -- (263.7993,539.1780)
  		(274.5585,542.6740) -- (263.7993,539.1780)
  		(281.2081,563.1396) -- (287.8577,572.2921)
  		(239.7401,572.2921) -- (239.7401,583.6052)
  		(239.7401,583.6052) -- (228.9809,587.1012)
  		(228.9809,587.1012) -- (218.2209,583.6052)
  		(224.8713,592.7577) -- (218.2209,583.6052)
  		(224.8713,592.7577) -- (218.2209,601.9103)
  		(253.0393,613.2234) -- (253.0393,601.9103)
  		(263.7993,605.4062) -- (253.0393,601.9103)
  		(253.0393,613.2234) -- (242.2801,616.7194)
  		(242.2801,616.7194) -- (235.6305,607.5668)
  		(281.2081,622.3760) -- (291.9673,625.8719)
  		(291.9673,596.2537) -- (298.6169,587.1012)
  		(309.3769,601.9103) -- (302.7265,592.7577)
  		(302.7265,592.7577) -- (291.9673,596.2537)
  		(281.2081,622.3760) -- (274.5585,613.2234)
  		(274.5585,613.2234) -- (263.7985,616.7194)
  		(302.7265,629.3678) -- (313.4865,625.8719)
  		(313.4865,625.8719) -- (320.1361,635.0244)
  		(320.1361,635.0244) -- (320.1361,646.3376)
  		(326.7857,637.1850) -- (320.1361,646.3376)
  		(326.7857,637.1850) -- (337.5449,640.6809)
  		(337.5449,604.0709) -- (326.7857,607.5668)
  		(337.5449,622.3760) -- (326.7857,625.8719)
  		(326.7857,625.8719) -- (337.5449,629.3678)
  		(320.1361,605.4062) -- (326.7857,596.2537)
  		(326.7857,596.2537) -- (326.7857,607.5668)
  		(337.5449,604.0709) -- (344.1945,613.2234)
  		(344.1945,613.2234) -- (337.5449,622.3760)
  		(337.5449,629.3678) -- (337.5449,640.6809)
  		(309.3769,601.9103) -- (320.1361,605.4062)
  		(263.7985,616.7194) -- (263.7993,605.4062)
  		(218.2209,601.9103) -- (228.9809,605.4062)
  		(228.9809,605.4062) -- (235.6305,596.2537)
  		(235.6305,596.2537) -- (235.6305,607.5668)
  		(263.7993,605.4062) -- (257.1489,596.2537)
  		(257.1489,596.2537) -- (263.7993,587.1012)
  		(302.7265,629.3679) -- (291.9673,625.8719)
  		(263.7993,605.4062) -- (274.5585,601.9103)
  		(281.2081,592.7577) -- (274.5585,601.9103);
  	\end{tikzpicture}
  	\caption{Arrangements for $T_1,T_2$ and $T_2,T_3$.}
  \end{figure}
  
  \begin{figure}[!ht]
  	\centering
  	\begin{tikzpicture}
  		[y=-1.0pt, x=1.0pt, scale=0.8]
  		\draw[red, dashed, line width=0.5pt]
  		(79.5906,537.2860) -- (5.8442,638.7884)
  		(5.8442,638.7884) -- (79.5906,662.7500)
  		(79.5906,662.7500) -- (153.3370,638.7884)
  		(153.3370,638.7884) -- (79.5906,537.2860)
  		(264.8473,537.5029) -- (338.5937,639.0053)
  		(338.5937,639.0053) -- (264.8473,615.0437)
  		(264.8473,615.0437) -- (191.1017,639.0053)
  		(264.8473,537.5029) -- (191.1017,639.0053);
  		
  		\draw[line width=0.5pt]
  		(90.3498,629.6364) -- (96.9994,638.7884)
  		(96.9994,638.7884) -- (86.2402,635.2924)
  		(96.9994,638.7884) -- (90.3498,647.9412)
  		(96.9994,638.7884) -- (103.6498,647.9412)
  		(96.9994,638.7884) -- (103.6498,647.9412)
  		(40.6626,609.1708) -- (51.4218,605.6748)
  		(29.9034,605.6748) -- (40.6626,602.1788)
  		(40.6626,602.1788) -- (51.4218,605.6748)
  		(40.6626,590.8652) -- (40.6626,602.1788)
  		(79.5906,651.4372) -- (90.3498,647.9412)
  		(79.5906,662.7500) -- (90.3498,659.2540)
  		(90.3498,659.2540) -- (90.3498,647.9412)
  		(90.3498,629.6364) -- (79.5906,626.1404)
  		(90.3498,629.6364) -- (90.3498,618.3228)
  		(90.3498,618.3228) -- (79.5906,614.8268)
  		(79.5906,626.1404) -- (79.5906,614.8268)
  		(68.8314,618.3228) -- (79.5906,614.8268)
  		(68.8314,618.3228) -- (58.0714,614.8268)
  		(51.4218,623.9796) -- (58.0714,614.8268)
  		(58.0714,614.8268) -- (51.4218,605.6748)
  		(44.7722,614.8268) -- (51.4218,605.6748)
  		(79.5906,644.4452) -- (90.3498,647.9412)
  		(103.6498,647.9412) -- (96.9994,657.0940)
  		(96.9994,638.7884) -- (90.3498,647.9412)
  		(96.9994,657.0940) -- (90.3498,647.9412)
  		(29.9034,605.6748) -- (29.9034,594.3612)
  		(40.6626,590.8652) -- (29.9034,594.3612)
  		(68.8314,618.3228) -- (62.1818,609.1708)
  		(62.1818,609.1708) -- (51.4218,605.6748)
  		(23.2538,644.4452) -- (34.0130,647.9412)
  		(34.0130,647.9412) -- (44.7722,651.4372)
  		(16.6042,623.9796) -- (23.2538,614.8268)
  		(23.2538,626.1404) -- (23.2538,614.8268)
  		(34.0130,647.9412) -- (27.3634,638.7884)
  		(27.3634,638.7884) -- (34.0130,629.6364)
  		(34.0130,629.6364) -- (23.2538,626.1404)
  		(34.0130,647.9412) -- (40.6626,638.7884)
  		(23.2538,614.8268) -- (34.0130,618.3228)
  		(34.0130,618.3228) -- (34.0130,629.6364)
  		(34.0130,618.3228) -- (40.6626,609.1708)
  		(40.6626,609.1708) -- (29.9034,605.6748)
  		(40.6626,627.4756) -- (40.6626,638.7884)
  		(34.0130,618.3228) -- (40.6626,627.4756)
  		(40.6626,627.4756) -- (51.4218,623.9796)
  		(51.4218,623.9796) -- (44.7722,614.8268)
  		(44.7722,614.8268) -- (34.0130,618.3228)
  		(51.4218,623.9796) -- (51.4218,635.2924)
  		(51.4218,635.2924) -- (40.6626,638.7884)
  		(51.4218,635.2924) -- (62.1818,638.7884)
  		(51.4218,642.2844) -- (62.1818,638.7884)
  		(23.2538,644.4452) -- (16.6042,635.2924)
  		(16.6042,635.2924) -- (16.6042,623.9796)
  		(16.6042,623.9796) -- (5.8442,627.4756)
  		(16.6042,635.2924) -- (5.8442,638.7884)
  		(5.8442,627.4756) -- (5.8442,638.7884)
  		(44.7722,651.4372) -- (51.4218,642.2844)
  		(51.4218,642.2844) -- (40.6626,638.7884)
  		(34.0130,629.6364) -- (40.6626,638.7884)
  		(34.0130,647.9412) -- (40.6626,638.7884)
  		(16.6042,635.2924) -- (23.2538,626.1404)
  		(16.6042,635.2924) -- (27.3634,638.7884)
  		(23.2538,614.8268) -- (29.9034,605.6748)
  		(44.7722,651.4372) -- (55.5322,647.9412)
  		(55.5322,647.9412) -- (62.1818,638.7884)
  		(55.5322,647.9412) -- (62.1818,657.0940)
  		(62.1818,657.0940) -- (62.1818,668.4068)
  		(62.1818,657.0940) -- (68.8314,647.9412)
  		(68.8314,647.9412) -- (62.1818,638.7884)
  		(68.8314,647.9412) -- (68.8314,659.2540)
  		(68.8314,659.2540) -- (62.1818,668.4068)
  		(68.8314,659.2540) -- (79.5906,662.7500)
  		(62.1818,627.4756) -- (62.1818,638.7884)
  		(68.8314,629.6364) -- (62.1818,638.7884)
  		(79.5906,626.1404) -- (68.8314,629.6364)
  		(79.5906,626.1404) -- (72.9410,635.2924)
  		(72.9410,635.2924) -- (62.1818,638.7884)
  		(72.9410,635.2924) -- (79.5906,644.4452)
  		(79.5906,644.4452) -- (68.8314,647.9412)
  		(68.8314,647.9412) -- (79.5906,651.4372)
  		(62.1818,627.4756) -- (68.8314,618.3228)
  		(68.8314,618.3228) -- (68.8314,629.6364)
  		(79.5906,626.1404) -- (86.2402,635.2924)
  		(86.2402,635.2924) -- (79.5906,644.4452)
  		(79.5906,651.4372) -- (79.5906,662.7500)
  		(51.4218,623.9796) -- (62.1818,627.4756)
  		(79.5906,596.5220) -- (79.5906,585.2092)
  		(79.5906,585.2092) -- (79.5906,573.8956)
  		(62.1818,602.1788) -- (51.4218,605.6748)
  		(79.5906,585.2092) -- (72.9410,594.3612)
  		(72.9410,594.3612) -- (62.1818,590.8652)
  		(62.1818,590.8652) -- (62.1818,602.1788)
  		(79.5906,585.2092) -- (68.8314,581.7132)
  		(51.4218,605.6748) -- (51.4218,594.3612)
  		(51.4218,594.3612) -- (62.1818,590.8652)
  		(51.4218,594.3612) -- (40.6626,590.8652)
  		(58.0722,585.2092) -- (68.8314,581.7132)
  		(51.4218,594.3612) -- (58.0722,585.2092)
  		(58.0722,585.2092) -- (51.4218,576.0564)
  		(51.4218,576.0564) -- (44.7722,585.2092)
  		(44.7722,585.2092) -- (51.4218,594.3612)
  		(51.4218,576.0564) -- (62.1818,572.5604)
  		(62.1818,572.5604) -- (68.8314,581.7132)
  		(62.1818,572.5604) -- (62.1818,561.2476)
  		(68.8314,570.3996) -- (62.1818,561.2476)
  		(79.5906,596.5220) -- (72.9410,605.6748)
  		(72.9410,605.6748) -- (62.1818,609.1708)
  		(72.9410,605.6748) -- (79.5906,614.8268)
  		(79.5906,573.8956) -- (68.8314,570.3996)
  		(68.8314,570.3996) -- (68.8314,581.7132)
  		(62.1818,590.8652) -- (68.8314,581.7132)
  		(79.5906,585.2092) -- (68.8314,581.7132)
  		(72.9410,605.6748) -- (62.1818,602.1788)
  		(72.9410,605.6748) -- (72.9410,594.3612)
  		(79.5906,573.8956) -- (72.9410,564.7436)
  		(72.9410,564.7436) -- (62.1818,561.2476)
  		(72.9410,564.7436) -- (79.5906,555.5908)
  		(79.5906,555.5908) -- (90.3498,552.0948)
  		(79.5906,555.5908) -- (68.8314,552.0948)
  		(68.8314,552.0948) -- (62.1818,561.2476)
  		(68.8314,552.0948) -- (79.5906,548.5988)
  		(79.5906,548.5988) -- (90.3498,552.0948)
  		(79.5906,548.5988) -- (79.5906,537.2860)
  		(51.4218,564.7436) -- (62.1818,561.2476)
  		(51.4218,557.7516) -- (62.1818,561.2476)
  		(44.7722,548.5988) -- (51.4218,557.7516)
  		(44.7722,548.5988) -- (55.5322,552.0948)
  		(55.5322,552.0948) -- (62.1818,561.2476)
  		(55.5322,552.0948) -- (62.1818,542.9420)
  		(62.1818,542.9420) -- (68.8314,552.0948)
  		(68.8314,552.0948) -- (68.8314,540.7820)
  		(51.4218,564.7436) -- (40.6626,561.2476)
  		(40.6626,561.2476) -- (51.4218,557.7516)
  		(44.7722,548.5988) -- (51.4218,539.4468)
  		(51.4218,539.4468) -- (62.1818,542.9420)
  		(68.8314,540.7820) -- (79.5906,537.2860)
  		(51.4218,576.0564) -- (51.4218,564.7436)
  		(79.5906,585.2092) -- (86.2402,594.3612)
  		(79.5906,585.2092) -- (90.3498,581.7132)
  		(107.7594,642.2844) -- (97.0002,638.7884)
  		(114.4090,651.4372) -- (103.6498,647.9412)
  		(101.1098,585.2092) -- (90.3498,581.7132)
  		(107.7594,576.0564) -- (97.0002,572.5604)
  		(97.0002,572.5604) -- (90.3498,581.7132)
  		(79.5906,596.5220) -- (86.2402,605.6748)
  		(86.2402,605.6748) -- (79.5906,614.8268)
  		(90.3498,618.3228) -- (97.0002,627.4756)
  		(107.7594,623.9796) -- (97.0002,627.4756)
  		(97.0002,627.4756) -- (97.0002,638.7884)
  		(107.7594,635.2924) -- (97.0002,638.7884)
  		(97.0002,590.8652) -- (90.3498,581.7132)
  		(79.5906,573.8956) -- (90.3498,570.3996)
  		(79.5906,585.2092) -- (90.3498,581.7132)
  		(90.3498,570.3996) -- (90.3498,581.7132)
  		(114.4090,651.4372) -- (107.7594,660.5900)
  		(97.0002,657.0940) -- (107.7594,660.5900)
  		(142.5778,623.9796) -- (135.9274,614.8268)
  		(135.9274,614.8268) -- (129.2778,605.6748)
  		(135.9274,644.4452) -- (125.1682,647.9412)
  		(131.8178,638.7884) -- (125.1682,647.9412)
  		(135.9274,614.8268) -- (135.9274,626.1404)
  		(135.9274,626.1404) -- (125.1682,629.6364)
  		(125.1682,629.6364) -- (131.8178,638.7884)
  		(135.9274,614.8268) -- (125.1682,618.3228)
  		(125.1682,647.9412) -- (118.5186,638.7884)
  		(118.5186,638.7884) -- (125.1682,629.6364)
  		(118.5186,638.7884) -- (107.7594,642.2844)
  		(107.7594,642.2844) -- (114.4090,651.4372)
  		(118.5186,627.4756) -- (125.1682,618.3228)
  		(118.5186,638.7884) -- (118.5186,627.4756)
  		(118.5186,627.4756) -- (107.7594,623.9796)
  		(107.7594,623.9796) -- (107.7594,635.2924)
  		(107.7594,635.2924) -- (118.5186,638.7884)
  		(107.7594,623.9796) -- (114.4090,614.8268)
  		(114.4090,614.8268) -- (125.1682,618.3228)
  		(114.4090,614.8268) -- (107.7594,605.6748)
  		(118.5186,609.1708) -- (107.7594,605.6748)
  		(142.5778,623.9796) -- (142.5778,635.2924)
  		(142.5778,635.2924) -- (135.9274,644.4452)
  		(135.9274,644.4452) -- (146.6874,647.9412)
  		(142.5778,635.2924) -- (153.3370,638.7884)
  		(146.6874,647.9412) -- (153.3370,638.7884)
  		(129.2778,605.6748) -- (118.5186,609.1708)
  		(118.5186,609.1708) -- (125.1682,618.3228)
  		(125.1682,629.6364) -- (125.1682,618.3228)
  		(135.9274,614.8268) -- (125.1682,618.3228)
  		(142.5778,635.2924) -- (131.8178,638.7884)
  		(142.5778,635.2924) -- (135.9274,626.1404)
  		(125.1682,647.9412) -- (114.4090,651.4372)
  		(129.2778,605.6748) -- (118.5186,602.1788)
  		(118.5186,602.1788) -- (107.7594,605.6748)
  		(118.5186,602.1788) -- (118.5186,590.8652)
  		(118.5186,590.8652) -- (125.1682,581.7132)
  		(118.5186,590.8652) -- (107.7594,594.3612)
  		(107.7594,594.3612) -- (107.7594,605.6748)
  		(107.7594,594.3612) -- (114.4090,585.2092)
  		(114.4090,585.2092) -- (125.1682,581.7132)
  		(114.4090,585.2092) -- (107.7594,576.0564)
  		(101.1098,614.8268) -- (107.7594,605.6748)
  		(96.9994,609.1708) -- (107.7594,605.6748)
  		(86.2402,605.6748) -- (96.9994,609.1708)
  		(86.2402,605.6748) -- (96.9994,602.1788)
  		(96.9994,602.1788) -- (107.7594,605.6748)
  		(96.9994,602.1788) -- (96.9994,590.8652)
  		(96.9994,590.8652) -- (107.7594,594.3612)
  		(107.7594,594.3612) -- (101.1098,585.2092)
  		(101.1098,614.8268) -- (90.3498,618.3228)
  		(90.3498,618.3228) -- (96.9994,609.1708)
  		(86.2402,605.6748) -- (86.2402,594.3612)
  		(86.2402,594.3612) -- (96.9994,590.8652)
  		(101.1098,585.2092) -- (107.7594,576.0564)
  		(107.7594,623.9796) -- (101.1098,614.8268)
  		(247.4385,620.7005) -- (236.6793,624.1965)
  		(236.6793,624.1965) -- (243.3289,615.0437)
  		(236.6793,624.1965) -- (230.0297,615.0437)
  		(236.6793,624.1965) -- (225.9193,627.6925)
  		(282.2569,579.7693) -- (282.2569,591.0821)
  		(288.9065,570.6165) -- (288.9065,581.9301)
  		(288.9065,581.9301) -- (282.2569,591.0821)
  		(299.6657,585.4261) -- (288.9065,581.9301)
  		(230.0297,603.7309) -- (230.0297,615.0437)
  		(219.2697,600.2349) -- (219.2697,611.5477)
  		(219.2697,611.5477) -- (230.0297,615.0437)
  		(247.4385,620.7005) -- (254.0881,611.5477)
  		(247.4385,620.7005) -- (258.1977,624.1965)
  		(258.1977,624.1965) -- (264.8473,615.0437)
  		(254.0881,611.5477) -- (264.8473,615.0437)
  		(264.8473,603.7309) -- (264.8473,615.0437)
  		(264.8473,603.7309) -- (271.4969,594.5781)
  		(264.8473,585.4261) -- (271.4969,594.5781)
  		(271.4969,594.5781) -- (282.2569,591.0821)
  		(275.6073,581.9301) -- (282.2569,591.0821)
  		(236.6793,605.8917) -- (230.0297,615.0437)
  		(225.9193,627.6925) -- (219.2697,618.5397)
  		(236.6793,624.1965) -- (230.0297,615.0437)
  		(219.2697,618.5397) -- (230.0297,615.0437)
  		(288.9065,570.6165) -- (299.6657,574.1125)
  		(299.6657,585.4261) -- (299.6657,574.1125)
  		(264.8473,603.7309) -- (275.6073,600.2349)
  		(275.6073,600.2349) -- (282.2569,591.0821)
  		(254.0881,552.3117) -- (247.4385,561.4645)
  		(247.4385,561.4645) -- (240.7889,570.6165)
  		(275.6073,552.3117) -- (282.2569,561.4645)
  		(271.4969,557.9685) -- (282.2569,561.4645)
  		(247.4385,561.4645) -- (258.1977,557.9685)
  		(258.1977,557.9685) -- (264.8473,567.1205)
  		(264.8473,567.1205) -- (271.4969,557.9685)
  		(247.4385,561.4645) -- (254.0881,570.6165)
  		(282.2569,561.4645) -- (275.6073,570.6165)
  		(275.6073,570.6165) -- (264.8473,567.1205)
  		(275.6073,570.6165) -- (282.2569,579.7693)
  		(282.2569,579.7693) -- (288.9065,570.6165)
  		(264.8473,574.1125) -- (254.0881,570.6165)
  		(275.6073,570.6165) -- (264.8473,574.1125)
  		(264.8473,574.1125) -- (264.8473,585.4261)
  		(264.8473,585.4261) -- (275.6073,581.9301)
  		(275.6073,581.9301) -- (275.6073,570.6165)
  		(264.8473,585.4261) -- (254.0881,581.9301)
  		(254.0881,581.9301) -- (254.0881,570.6165)
  		(254.0881,581.9301) -- (247.4385,591.0821)
  		(247.4385,579.7693) -- (247.4385,591.0821)
  		(254.0881,552.3117) -- (264.8473,548.8157)
  		(264.8473,548.8157) -- (275.6073,552.3117)
  		(275.6073,552.3117) -- (275.6073,540.9989)
  		(264.8473,548.8157) -- (264.8473,537.5029)
  		(275.6073,540.9989) -- (264.8473,537.5029)
  		(240.7889,570.6165) -- (247.4385,579.7693)
  		(247.4385,579.7693) -- (254.0881,570.6165)
  		(264.8473,567.1205) -- (254.0881,570.6165)
  		(247.4385,561.4645) -- (254.0881,570.6165)
  		(264.8473,548.8157) -- (271.4969,557.9685)
  		(264.8473,548.8157) -- (258.1977,557.9685)
  		(282.2569,561.4645) -- (288.9065,570.6165)
  		(240.7889,570.6165) -- (240.7889,581.9301)
  		(240.7889,581.9301) -- (247.4385,591.0821)
  		(240.7889,581.9301) -- (230.0297,585.4261)
  		(230.0297,585.4261) -- (219.2697,581.9301)
  		(230.0297,585.4261) -- (236.6793,594.5781)
  		(236.6793,594.5781) -- (247.4385,591.0821)
  		(236.6793,594.5781) -- (225.9193,591.0821)
  		(225.9193,591.0821) -- (219.2697,581.9301)
  		(225.9193,591.0821) -- (219.2697,600.2349)
  		(258.1977,594.5781) -- (247.4385,591.0821)
  		(254.0881,600.2349) -- (247.4385,591.0821)
  		(254.0881,611.5477) -- (254.0881,600.2349)
  		(254.0881,611.5477) -- (247.4385,602.3957)
  		(247.4385,602.3957) -- (247.4385,591.0821)
  		(247.4385,602.3957) -- (236.6793,605.8917)
  		(236.6793,605.8917) -- (236.6793,594.5781)
  		(236.6793,594.5781) -- (230.0297,603.7309)
  		(258.1977,594.5781) -- (264.8473,603.7309)
  		(264.8473,603.7309) -- (254.0881,600.2349)
  		(254.0881,611.5477) -- (243.3289,615.0437)
  		(243.3289,615.0437) -- (236.6793,605.8917)
  		(230.0297,603.7309) -- (219.2697,600.2349)
  		(264.8473,585.4261) -- (258.1977,594.5781)
  		(282.2569,620.7005) -- (293.0161,624.1965)
  		(293.0161,624.1965) -- (303.7753,627.6925)
  		(282.2569,602.3957) -- (282.2569,591.0821)
  		(293.0161,624.1965) -- (286.3665,615.0437)
  		(286.3665,615.0437) -- (293.0161,605.8917)
  		(293.0161,605.8917) -- (282.2569,602.3957)
  		(293.0161,624.1965) -- (299.6657,615.0437)
  		(282.2569,591.0821) -- (293.0161,594.5781)
  		(293.0161,594.5781) -- (293.0161,605.8917)
  		(293.0161,594.5781) -- (299.6657,585.4261)
  		(299.6657,603.7309) -- (299.6657,615.0437)
  		(293.0161,594.5781) -- (299.6657,603.7309)
  		(299.6657,603.7309) -- (310.4249,600.2349)
  		(310.4249,600.2349) -- (303.7753,591.0821)
  		(303.7753,591.0821) -- (293.0161,594.5781)
  		(310.4249,600.2349) -- (310.4249,611.5477)
  		(310.4249,611.5477) -- (299.6657,615.0437)
  		(310.4249,611.5477) -- (321.1849,615.0437)
  		(310.4249,618.5397) -- (321.1849,615.0437)
  		(282.2569,620.7005) -- (275.6073,611.5477)
  		(275.6073,611.5477) -- (275.6073,600.2349)
  		(275.6073,611.5477) -- (264.8473,615.0437)
  		(303.7753,627.6925) -- (310.4249,618.5397)
  		(310.4249,618.5397) -- (299.6657,615.0437)
  		(293.0161,605.8917) -- (299.6657,615.0437)
  		(293.0161,624.1965) -- (299.6657,615.0437)
  		(275.6073,611.5477) -- (282.2569,602.3957)
  		(275.6073,611.5477) -- (286.3665,615.0437)
  		(303.7753,627.6925) -- (314.5353,624.1965)
  		(314.5353,624.1965) -- (321.1849,615.0437)
  		(314.5353,624.1965) -- (321.1849,633.3493)
  		(321.1849,633.3493) -- (321.1849,644.6621)
  		(321.1849,633.3493) -- (327.8345,624.1965)
  		(327.8345,624.1965) -- (321.1849,615.0437)
  		(327.8345,624.1965) -- (327.8345,635.5093)
  		(327.8345,635.5093) -- (321.1849,644.6621)
  		(327.8345,635.5093) -- (338.5937,639.0053)
  		(321.1849,603.7309) -- (321.1849,615.0437)
  		(327.8345,605.8917) -- (321.1849,615.0437)
  		(338.5937,602.3957) -- (327.8345,605.8917)
  		(338.5937,602.3957) -- (331.9441,611.5477)
  		(331.9441,611.5477) -- (321.1849,615.0437)
  		(331.9441,611.5477) -- (338.5937,620.7005)
  		(338.5937,620.7005) -- (327.8345,624.1965)
  		(327.8345,624.1965) -- (338.5937,627.6925)
  		(321.1849,603.7309) -- (327.8345,594.5781)
  		(327.8345,594.5781) -- (327.8345,605.8917)
  		(338.5937,602.3957) -- (345.2433,611.5477)
  		(345.2433,611.5477) -- (338.5937,620.7005)
  		(338.5937,627.6925) -- (338.5937,639.0053)
  		(310.4249,600.2349) -- (321.1849,603.7309);
  	\end{tikzpicture}
  	\caption{Tiles of Figure 2 with decoration.}
  \end{figure}
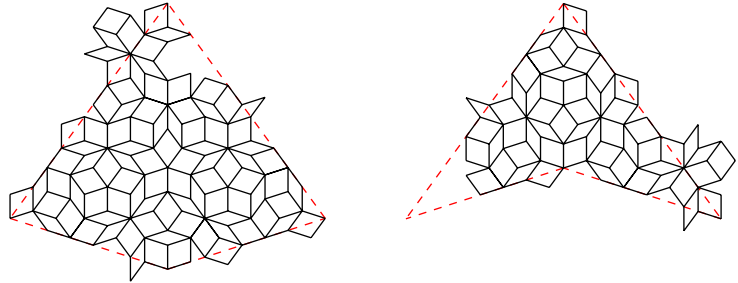
  
  The dashed kite and dart shapes in Figure 3 can be completely decorated or filled by rhombs, including half rhombs with their diagonals on the dashed edges. The arrangements in Figures 2 and 3 can also be interpreted as modifications of the kite and dart shape, as mentioned in \cite{Gruenbaum} on page 539, where the long and short sides of kite and dart are replaced by two J-curves.

  \section{Variants of the Prototiles} 
  
  From $T_1$ and $T_2$ we can get two new undecorated prototiles, nicknamed Snake and Dog, by shifting two outside half thin rhombs to the corresponding gaps (Fig. 4)\footnote{The black dots within the tiles should represent the animals' eyes only.}.
  
  \begin{figure}[!ht]
  	\centering
  	\begin{tikzpicture}
  		[y=1.0pt, x=1.0pt, yscale=-0.9, xscale=0.9]
  		\draw[line width=0.5pt]
  		(11.9980,602.1048) -- (15.4940,591.3457)
  		(18.9900,602.1048) -- (30.3028,602.1048)
  		(15.4940,591.3457) -- (18.9900,580.5863)
  		(30.3028,636.9232) -- (26.8068,626.1640)
  		(39.4556,608.7544) -- (30.3028,602.1048)
  		(24.6460,654.3328) -- (28.1420,665.0920)
  		(24.6460,654.3328) -- (35.9596,654.3328)
  		(18.9900,658.4424) -- (15.4940,647.6824)
  		(9.8372,665.0920) -- (18.9900,658.4424)
  		(39.4556,608.7544) -- (35.9596,597.9954)
  		(35.9596,597.9954) -- (45.1124,591.3457)
  		(45.1124,591.3457) -- (35.9596,584.6960)
  		(35.9596,584.6960) -- (26.8068,591.3457)
  		(11.9980,602.1048) -- (2.8452,608.7544)
  		(2.8452,608.7544) -- (6.3412,619.5144)
  		(24.6460,619.5144) -- (15.4940,626.1640)
  		(15.4940,612.8648) -- (24.6460,619.5144)
  		(15.4940,626.1640) -- (6.3412,619.5144)
  		(15.4940,626.1640) -- (26.8068,626.1640)
  		(35.9596,654.3328) -- (26.8068,647.6824)
  		(30.3028,636.9232) -- (26.8068,647.6824)
  		(15.4940,612.8648) -- (18.9900,602.1048)
  		(18.9900,580.5863) -- (30.3028,580.5863)
  		(30.3028,580.5863) -- (26.8068,591.3457)
  		(28.1420,665.0920) -- (18.9900,671.7416)
  		(9.8372,665.0920) -- (18.9900,671.7416)
  		(15.4940,626.1640) -- (11.9980,636.9232)
  		(11.9980,636.9232) -- (15.4940,647.6824)
  		(168.5529,647.7092) -- (165.0569,636.9500)
  		(165.0569,636.9500) -- (161.5609,626.1908)
  		(155.9041,665.1180) -- (144.5913,665.1180)
  		(141.0953,654.3588) -- (129.7825,654.3588)
  		(129.7825,654.3588) -- (133.2785,665.1180)
  		(135.4385,636.9500) -- (131.9425,647.7092)
  		(131.9425,647.7092) -- (141.0953,654.3588)
  		(168.5529,647.7092) -- (165.0569,658.4684)
  		(155.9041,665.1180) -- (165.0569,671.7684)
  		(165.0569,658.4684) -- (174.2097,665.1180)
  		(165.0569,671.7684) -- (174.2097,665.1180)
  		(144.5913,665.1180) -- (133.2785,665.1180)
  		(161.5609,626.1908) -- (152.4081,619.5404)
  		(152.4081,619.5404) -- (155.9041,608.7812)
  		(155.9041,608.7812) -- (165.0569,602.1316)
  		(153.7441,602.1316) -- (165.0569,602.1316)
  		(153.7441,602.1316) -- (150.2481,591.3722)
  		(120.6297,612.8908) -- (129.7825,619.5404)
  		(135.4385,602.1316) -- (144.5913,608.7812)
  		(144.5913,608.7812) -- (141.0953,598.0219)
  		(131.9425,626.1908) -- (120.6297,626.1908)
  		(120.6297,626.1908) -- (129.7825,619.5404)
  		(120.6297,612.8908) -- (124.1257,602.1316)
  		(124.1257,602.1316) -- (135.4385,602.1316)
  		(141.0953,598.0219) -- (150.2481,591.3722)
  		(135.4385,636.9500) -- (131.9425,626.1908)
  		(68.6625,602.2838) -- (72.1585,591.5244)
  		(75.6545,602.2838) -- (86.9681,602.2838)
  		(72.1585,591.5244) -- (75.6545,580.7651)
  		(86.9681,637.1022) -- (83.4721,626.3430)
  		(81.3113,654.5110) -- (84.8073,665.2710)
  		(75.6545,658.6206) -- (72.1585,647.8614)
  		(66.5025,665.2710) -- (75.6545,658.6206)
  		(92.6241,598.1742) -- (101.7769,591.5244)
  		(101.7769,591.5244) -- (92.6241,584.8748)
  		(68.6625,602.2838) -- (59.5105,608.9334)
  		(59.5105,608.9334) -- (63.0065,619.6934)
  		(72.1585,613.0430) -- (81.3113,619.6934)
  		(72.1585,626.3430) -- (63.0065,619.6934)
  		(86.9681,637.1022) -- (83.4721,647.8614)
  		(72.1585,613.0430) -- (75.6545,602.2838)
  		(75.6545,580.7651) -- (86.9681,580.7651)
  		(84.8073,665.2710) -- (75.6545,671.9206)
  		(66.5025,665.2710) -- (75.6545,671.9206)
  		(72.1585,626.3430) -- (68.6625,637.1022)
  		(68.6625,637.1022) -- (72.1585,647.8614)
  		(86.9681,580.7651) -- (92.6241,584.8748)
  		(81.3113,619.6934) -- (83.4721,626.3430)
  		(81.3113,654.5110) -- (83.4721,647.8614)
  		(86.9681,602.2838) -- (92.6241,598.1742)
  		(231.1647,648.2160) -- (227.6687,637.4560)
  		(227.6687,637.4560) -- (224.1727,626.6968)
  		(218.5159,665.6248) -- (207.2031,665.6248)
  		(192.3935,654.8656) -- (195.8895,665.6248)
  		(198.0503,637.4560) -- (194.5543,648.2160)
  		(231.1647,648.2160) -- (227.6687,658.9752)
  		(218.5159,665.6248) -- (227.6687,672.2744)
  		(227.6687,658.9752) -- (236.8215,665.6248)
  		(227.6687,672.2744) -- (236.8215,665.6248)
  		(207.2031,665.6248) -- (195.8895,665.6248)
  		(224.1727,626.6968) -- (215.0199,620.0472)
  		(215.0199,620.0472) -- (218.5159,609.2880)
  		(216.3551,602.6381) -- (212.8599,591.8786)
  		(183.2415,613.3976) -- (192.3935,620.0472)
  		(183.2415,613.3976) -- (186.7375,602.6381)
  		(186.7375,602.6381) -- (198.0503,602.6381)
  		(203.7071,598.5284) -- (212.8599,591.8786)
  		(198.0503,637.4560) -- (194.5543,626.6968)
  		(192.3935,654.8656) -- (194.5543,648.2160)
  		(198.0503,602.6381) -- (203.7071,598.5284)
  		(218.5159,609.2880) -- (216.3551,602.6381)
  		(194.5543,626.6968) -- (192.3935,620.0472);
  		
  		\fill[line width=0.1pt]  (87,588) circle[radius=0.05cm];
  		\fill[line width=0.1pt]  (87,595) circle[radius=0.05cm];
  		\fill[line width=0.1pt] (203,606) circle[radius=0.05cm];
  		
  		\draw[red, <-] (40,580) to [bend right] (50,610);
  		\draw[red, <-] (30,623) to [bend right] (40,655);
  		\draw[red, <-] (140,590) to [bend right] (165,595);
  		\draw[red, ->] (118,630) to [bend left] (127,650);
  		
  		\small
  		\coordinate[label=right:$T_1$] (a) at (13,685);
  		\coordinate[label=right:Snake] (a) at (62,685);
  		\coordinate[label=right:$T_2$] (a) at (142,685);
  		\coordinate[label=right:Dog] (a) at (200,685);
  		\normalsize
  		
  		\coordinate[label=right:$\rightarrow$] (a) at (42,630);
  		\coordinate[label=right:$\rightarrow$] (a) at (168,630);
  	\end{tikzpicture}	
  	\caption{Building the snake and dog tiles.}
  \end{figure}
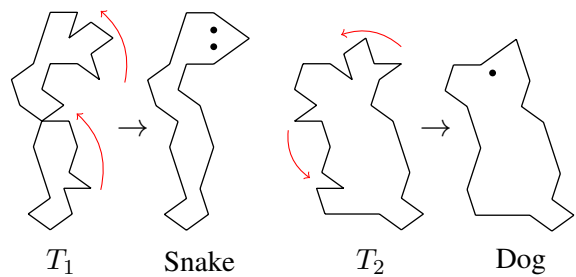
  
  Again the snake and dog tiles can be transformed into further and smoother shapes given by a hexagon and a pentagon, both concave and irregular. Figures 5 and 6 show this process from the left to the right by shifting three or four half rhombs as shown. $T_8$ and $T_9$, as well as Snake and Dog, are balanced tiles with the same acreage as $T_1$ and $T_2$.
  
  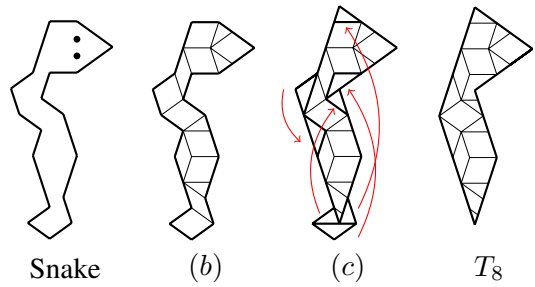
\begin{figure}[!ht]
  	\centering
  	\begin{tikzpicture}
  		[y=1.0pt, x=1.0pt, yscale=-0.9, xscale=0.9]
  		\draw[line width=0.01pt]
  		(126.8654,666.3076) -- (136.0182,659.6578)
  		(126.8654,666.3076) -- (136.0182,672.9572)
  		(132.5222,592.5613) -- (136.0182,603.3208)
  		(132.5222,592.5613) -- (143.8358,592.5613)
  		(141.6750,655.5481) -- (132.5222,648.8984)
  		(152.9878,599.2110) -- (143.8358,592.5613)
  		(152.9878,585.9116) -- (143.8358,592.5613)
  		(129.0262,603.3208) -- (132.5222,614.0801)
  		(141.6750,620.7299) -- (132.5222,627.3796)
  		(132.5222,614.0801) -- (123.3702,620.7299)
  		(132.5222,627.3796) -- (143.8358,627.3796)
  		(132.5222,627.3796) -- (136.0182,638.1390)
  		(147.3318,638.1390) -- (136.0182,638.1390)
  		(136.0182,638.1390) -- (132.5222,648.8984)
  		(143.8358,648.8984) -- (132.5222,648.8984)
  		(147.3318,603.3208) -- (143.8358,592.5613)
  		(147.3318,581.8019) -- (143.8358,592.5613)
  		(136.0182,659.6578) -- (139.5142,662.1978)
  		(129.0262,603.3208) -- (119.8742,609.9704)
  		(69.6954,603.3858) -- (73.1914,592.6264)
  		(73.1914,592.6264) -- (76.6866,603.3858)
  		(76.6866,603.3858) -- (88.0002,603.3858)
  		(73.1914,592.6264) -- (84.5042,592.6264)
  		(73.1914,592.6264) -- (76.6866,581.8669)
  		(88.0002,638.2040) -- (84.5042,627.4446)
  		(82.3434,655.6131) -- (85.8394,666.3726)
  		(82.3434,655.6131) -- (73.1914,648.9634)
  		(85.8394,666.3726) -- (76.6866,659.7228)
  		(76.6866,659.7228) -- (73.1914,648.9634)
  		(67.5346,666.3726) -- (76.6866,659.7228)
  		(93.6570,599.2760) -- (84.5042,592.6264)
  		(93.6570,599.2760) -- (102.8090,592.6264)
  		(102.8090,592.6264) -- (93.6570,585.9767)
  		(93.6570,585.9767) -- (84.5042,592.6264)
  		(69.6954,603.3858) -- (73.1914,614.1452)
  		(69.6954,603.3858) -- (60.5426,610.0355)
  		(60.5426,610.0355) -- (64.0386,620.7948)
  		(82.3434,620.7948) -- (73.1914,627.4446)
  		(73.1914,614.1452) -- (82.3434,620.7948)
  		(73.1914,614.1452) -- (64.0386,620.7948)
  		(73.1914,627.4446) -- (64.0386,620.7948)
  		(73.1914,627.4446) -- (84.5042,627.4446)
  		(73.1914,627.4446) -- (76.6866,638.2040)
  		(88.0002,638.2040) -- (76.6866,638.2040)
  		(76.6866,638.2040) -- (73.1914,648.9634)
  		(84.5042,648.9634) -- (73.1914,648.9634)
  		(88.0002,638.2040) -- (84.5042,648.9634)
  		(88.0002,603.3858) -- (84.5042,592.6264)
  		(73.1914,614.1452) -- (76.6866,603.3858)
  		(76.6866,581.8669) -- (88.0002,581.8669)
  		(88.0002,581.8669) -- (84.5042,592.6264)
  		(85.8394,666.3726) -- (76.6866,673.0223)
  		(67.5346,666.3726) -- (76.6866,673.0223)
  		(73.1914,627.4446) -- (69.6954,638.2040)
  		(69.6954,638.2040) -- (73.1914,648.9634)
  		(88.0002,581.8669) -- (93.6570,585.9767)
  		(82.3434,620.7948) -- (84.5042,627.4446)
  		(82.3434,655.6131) -- (84.5042,648.9634)
  		(88.0002,603.3858) -- (93.6570,599.2760)
  		(138.1790,666.3076) -- (145.1710,666.3076)
  		(119.8742,609.9704) -- (129.0262,603.3208)
  		(123.3702,620.7299) -- (132.5222,627.3796)
  		(188.8496,593.0223) -- (192.3456,603.7817)
  		(192.3456,603.7817) -- (203.6592,603.7817)
  		(188.8496,593.0223) -- (200.1632,593.0223)
  		(198.0024,656.0091) -- (188.8496,649.3594)
  		(192.3456,660.1188) -- (188.8496,649.3594)
  		(209.3152,599.6720) -- (200.1632,593.0223)
  		(209.3152,599.6720) -- (218.4680,593.0223)
  		(209.3152,599.6720) -- (218.4680,593.0223)
  		(209.3152,586.3726) -- (200.1632,593.0223)
  		(185.3536,603.7817) -- (188.8496,614.5411)
  		(198.0024,621.1908) -- (188.8496,627.8406)
  		(188.8496,614.5411) -- (198.0024,621.1908)
  		(188.8496,614.5411) -- (179.6976,621.1908)
  		(188.8496,627.8406) -- (179.6976,621.1908)
  		(188.8496,627.8406) -- (185.3536,638.6000)
  		(188.8496,627.8406) -- (200.1632,627.8406)
  		(188.8496,627.8406) -- (192.3456,638.6000)
  		(203.6592,638.6000) -- (192.3456,638.6000)
  		(192.3456,638.6000) -- (188.8496,649.3594)
  		(200.1632,649.3594) -- (188.8496,649.3594)
  		(203.6592,603.7817) -- (200.1632,593.0223)
  		(188.8496,614.5411) -- (192.3456,603.7817)
  		(192.3456,582.2629) -- (203.6592,582.2629)
  		(203.6592,582.2629) -- (200.1632,593.0223)
  		(188.8496,614.5411) -- (195.8416,614.5411)
  		(192.3456,660.1188) -- (195.8416,662.6588)
  		(194.5064,610.4315) -- (218.4680,593.0223)
  		(194.5064,575.6135) -- (218.4680,593.0223)
  		(179.6976,621.1908) -- (194.5064,666.7685)
  		(194.5064,610.4315) -- (203.6592,638.6000)
  		(179.6976,621.1908) -- (194.5064,575.6135)
  		(194.5064,666.7685) -- (203.6592,638.6000)
  		(138.1790,575.1518) -- (136.0182,581.8019)
  		(136.0182,581.8019) -- (129.0262,603.3208)
  		(119.8742,609.9704) -- (123.3702,620.7299)
  		(123.3702,620.7299) -- (132.5222,648.8984)
  		(132.5222,648.8984) -- (136.0182,659.6578)
  		(132.5222,627.3796) -- (129.0262,638.1390)
  		(129.0262,603.3208) -- (123.3702,620.7299)
  		(147.3318,581.8019) -- (136.0182,581.8019)
  		(138.1790,575.1518) -- (147.3318,581.8019)
  		(147.3318,581.8019) -- (162.1406,592.5613)
  		(147.3318,638.1390) -- (138.1790,609.9704)
  		(138.1790,609.9704) -- (162.1406,592.5613)
  		(136.0182,603.3208) -- (147.3318,603.3208)
  		(141.6750,620.7299) -- (132.5222,614.0801)
  		(136.0182,603.3208) -- (132.5222,614.0801)
  		(132.5222,614.0801) -- (138.1790,609.9704)
  		(147.3318,638.1390) -- (139.5142,662.1978)
  		(139.5142,662.1978) -- (138.1790,666.3076)
  		(136.0182,672.9572) -- (145.1710,666.3076)
  		(145.1710,666.3076) -- (141.6750,655.5481)
  		(126.8654,666.3076) -- (138.1790,666.3076)
  		(136.0182,659.6578) -- (138.1790,666.3076);
		
  		\draw[line width=0.8pt]
  		(126.8654,666.3076) -- (136.0182,659.6578)
  		(126.8654,666.3076) -- (136.0182,672.9572)
  		(129.0262,603.3208) -- (119.8742,609.9704)
  		(69.6954,603.3858) -- (73.1914,592.6264)
  		(76.6866,603.3858) -- (88.0002,603.3858)
  		(73.1914,592.6264) -- (76.6866,581.8669)
  		(88.0002,638.2040) -- (84.5042,627.4446)
  		(82.3434,655.6131) -- (85.8394,666.3726)
  		(76.6866,659.7228) -- (73.1914,648.9634)
  		(67.5346,666.3726) -- (76.6866,659.7228)
  		(93.6570,599.2760) -- (102.8090,592.6264)
  		(102.8090,592.6264) -- (93.6570,585.9767)
  		(69.6954,603.3858) -- (60.5426,610.0355)
  		(60.5426,610.0355) -- (64.0386,620.7948)
  		(73.1914,614.1452) -- (82.3434,620.7948)
  		(73.1914,627.4446) -- (64.0386,620.7948)
  		(88.0002,638.2040) -- (84.5042,648.9634)
  		(73.1914,614.1452) -- (76.6866,603.3858)
  		(76.6866,581.8669) -- (88.0002,581.8669)
  		(85.8394,666.3726) -- (76.6866,673.0223)
  		(67.5346,666.3726) -- (76.6866,673.0223)
  		(73.1914,627.4446) -- (69.6954,638.2040)
  		(69.6954,638.2040) -- (73.1914,648.9634)
  		(88.0002,581.8669) -- (93.6570,585.9767)
  		(82.3434,620.7948) -- (84.5042,627.4446)
  		(82.3434,655.6131) -- (84.5042,648.9634)
  		(88.0002,603.3858) -- (93.6570,599.2760)
  		(10.6895,603.0017) -- (14.1855,592.2423)
  		(17.6815,603.0017) -- (28.9951,603.0017)
  		(14.1855,592.2423) -- (17.6815,581.4829)
  		(28.9951,637.8200) -- (25.4991,627.0606)
  		(23.3383,655.2291) -- (26.8343,665.9885)
  		(17.6815,659.3388) -- (14.1855,648.5794)
  		(8.5295,665.9885) -- (17.6815,659.3388)
  		(34.6511,598.8920) -- (43.8039,592.2423)
  		(43.8039,592.2423) -- (34.6511,585.5927)
  		(10.6895,603.0017) -- (1.5375,609.6514)
  		(1.5375,609.6514) -- (5.0335,620.4108)
  		(14.1855,613.7612) -- (23.3383,620.4108)
  		(14.1855,627.0606) -- (5.0335,620.4108)
  		(28.9951,637.8200) -- (25.4991,648.5794)
  		(14.1855,613.7612) -- (17.6815,603.0017)
  		(17.6815,581.4829) -- (28.9951,581.4829)
  		(26.8343,665.9885) -- (17.6815,672.6382)
  		(8.5295,665.9885) -- (17.6815,672.6382)
  		(14.1855,627.0606) -- (10.6895,637.8200)
  		(10.6895,637.8200) -- (14.1855,648.5794)
  		(28.9951,581.4829) -- (34.6511,585.5927)
  		(23.3383,620.4108) -- (25.4991,627.0606)
  		(23.3383,655.2291) -- (25.4991,648.5794)
  		(28.9951,603.0017) -- (34.6511,598.8920)
  		(138.1790,666.3076) -- (145.1710,666.3076)
  		(119.8742,609.9704) -- (129.0262,603.3208)
  		(123.3702,620.7299) -- (132.5222,627.3796)
  		(209.3152,599.6720) -- (218.4680,593.0223)
  		(209.3152,599.6720) -- (218.4680,593.0223)
  		(194.5064,610.4315) -- (218.4680,593.0223)
  		(194.5064,575.6135) -- (218.4680,593.0223)
  		(179.6976,621.1908) -- (194.5064,666.7685)
  		(194.5064,610.4315) -- (203.6592,638.6000)
  		(179.6976,621.1908) -- (194.5064,575.6135)
  		(194.5064,666.7685) -- (203.6592,638.6000)
  		(138.1790,575.1518) -- (136.0182,581.8019)
  		(136.0182,581.8019) -- (129.0262,603.3208)
  		(119.8742,609.9704) -- (123.3702,620.7299)
  		(123.3702,620.7299) -- (132.5222,648.8984)
  		(132.5222,648.8984) -- (136.0182,659.6578)
  		(132.5222,627.3796) -- (129.0262,638.1390)
  		(129.0262,603.3208) -- (123.3702,620.7299)
  		(147.3318,581.8019) -- (136.0182,581.8019)
  		(138.1790,575.1518) -- (147.3318,581.8019)
  		(147.3318,581.8019) -- (162.1406,592.5613)
  		(147.3318,638.1390) -- (138.1790,609.9704)
  		(138.1790,609.9704) -- (162.1406,592.5613)
  		(136.0182,603.3208) -- (147.3318,603.3208)
  		(141.6750,620.7299) -- (132.5222,614.0801)
  		(136.0182,603.3208) -- (132.5222,614.0801)
  		(132.5222,614.0801) -- (138.1790,609.9704)
  		(147.3318,638.1390) -- (139.5142,662.1978)
  		(139.5142,662.1978) -- (138.1790,666.3076)
  		(136.0182,672.9572) -- (145.1710,666.3076)
  		(145.1710,666.3076) -- (141.6750,655.5481)
  		(126.8654,666.3076) -- (138.1790,666.3076)
  		(136.0182,659.6578) -- (138.1790,666.3076);
  		
  		\fill[line width=0.1pt] (29,589) circle[radius=0.05cm];
  		\fill[line width=0.1pt] (29,596) circle[radius=0.05cm];
  		
  		\draw[red, <-] (141,584) to [bend right] (146,660);
  		\draw[red, <-] (142,610) to [bend right] (146,672);
  		\draw[red, ->] (115,610) to [bend left] (122,632);
  		\draw[red, <-] (136,618) to [bend left] (130,662);
  		
  		\small
  		\coordinate[label=right:Snake] (a) at (5,685);
  		\coordinate[label=right:$(b)$] (a) at (71,685);
  		\coordinate[label=right:$(c)$] (a) at (130,685);
  		\coordinate[label=right:$T_8$] (a) at (190,685);
  		\normalsize
  	\end{tikzpicture}	
  	\caption{Building $T_8$ from the Snake.}
  \end{figure}
  
  \begin{figure}[!ht]
  	\centering
  	\begin{tikzpicture}
  		[y=1.0pt, x=1.0pt, yscale=-0.9, xscale=0.9]
  		\draw[line width=0.01pt]
  		(112.6818,648.9364) -- (109.1858,638.1770)
  		(109.1858,638.1770) -- (105.6898,627.4175)
  		(100.0330,666.3455) -- (88.7202,666.3455)
  		(97.8722,659.6958) -- (88.7202,666.3455)
  		(109.1858,638.1770) -- (105.6898,648.9364)
  		(105.6898,648.9364) -- (94.3762,648.9364)
  		(94.3762,648.9364) -- (97.8722,659.6958)
  		(109.1858,638.1770) -- (97.8722,638.1770)
  		(88.7202,666.3455) -- (85.2242,655.5861)
  		(85.2242,655.5861) -- (94.3762,648.9364)
  		(85.2242,655.5861) -- (73.9106,655.5861)
  		(73.9106,655.5861) -- (77.4066,666.3455)
  		(88.7202,644.8266) -- (97.8722,638.1770)
  		(85.2242,655.5861) -- (88.7202,644.8266)
  		(88.7202,644.8266) -- (79.5674,638.1770)
  		(79.5674,638.1770) -- (76.0714,648.9364)
  		(76.0714,648.9364) -- (85.2242,655.5861)
  		(79.5674,638.1770) -- (88.7202,631.5272)
  		(88.7202,631.5272) -- (97.8722,638.1770)
  		(88.7202,631.5272) -- (85.2242,620.7678)
  		(94.3762,627.4175) -- (85.2242,620.7678)
  		(112.6818,648.9364) -- (109.1858,659.6958)
  		(109.1858,659.6958) -- (100.0330,666.3455)
  		(100.0330,666.3455) -- (109.1858,672.9952)
  		(109.1858,659.6958) -- (118.3378,666.3455)
  		(109.1858,672.9952) -- (118.3378,666.3455)
  		(105.6898,627.4175) -- (94.3762,627.4175)
  		(94.3762,627.4175) -- (97.8722,638.1770)
  		(94.3762,648.9364) -- (97.8722,638.1770)
  		(109.1858,659.6958) -- (97.8722,659.6958)
  		(109.1858,659.6958) -- (105.6898,648.9364)
  		(85.2242,655.5861) -- (94.3762,648.9364)
  		(88.7202,666.3455) -- (77.4066,666.3455)
  		(105.6898,627.4175) -- (96.5370,620.7678)
  		(96.5370,620.7678) -- (85.2242,620.7678)
  		(96.5370,620.7678) -- (100.0330,610.0084)
  		(100.0330,610.0084) -- (88.7202,610.0084)
  		(88.7202,610.0084) -- (85.2242,620.7678)
  		(88.7202,610.0084) -- (97.8722,603.3587)
  		(97.8722,603.3587) -- (94.3762,592.5993)
  		(76.0714,627.4175) -- (85.2242,620.7678)
  		(73.9106,620.7678) -- (85.2242,620.7678)
  		(64.7586,614.1182) -- (73.9106,620.7678)
  		(64.7586,614.1182) -- (76.0714,614.1182)
  		(76.0714,614.1182) -- (85.2242,620.7678)
  		(76.0714,614.1182) -- (79.5674,603.3587)
  		(79.5674,603.3587) -- (88.7202,610.0084)
  		(88.7202,610.0084) -- (85.2242,599.2490)
  		(64.7586,614.1182) -- (68.2546,603.3587)
  		(68.2546,603.3587) -- (79.5674,603.3587)
  		(85.2242,599.2490) -- (94.3762,592.5993)
  		(79.5674,638.1770) -- (76.0714,627.4175)
  		(73.9106,620.7678) -- (76.0714,627.4175)
  		(73.9106,655.5861) -- (76.0714,648.9364)
  		(79.5674,603.3587) -- (85.2242,599.2490)
  		(97.8722,603.3587) -- (100.0330,610.0084)
  		(163.0492,659.6893) -- (153.8964,666.3390)
  		(174.3628,638.1704) -- (170.8668,648.9298)
  		(170.8668,648.9298) -- (159.5532,648.9298)
  		(159.5532,648.9298) -- (163.0492,659.6893)
  		(174.3628,638.1704) -- (163.0492,638.1704)
  		(153.8964,666.3390) -- (150.4012,655.5795)
  		(150.4012,655.5795) -- (159.5532,648.9298)
  		(150.4012,655.5795) -- (139.0876,655.5795)
  		(153.8964,644.8201) -- (163.0492,638.1704)
  		(150.4012,655.5795) -- (153.8964,644.8201)
  		(153.8964,644.8201) -- (144.7444,638.1704)
  		(141.2484,648.9298) -- (150.4012,655.5795)
  		(144.7444,638.1704) -- (153.8964,631.5207)
  		(153.8964,631.5207) -- (163.0492,638.1704)
  		(153.8964,631.5207) -- (150.4012,620.7613)
  		(159.5532,627.4110) -- (150.4012,620.7613)
  		(177.8580,648.9298) -- (174.3628,659.6893)
  		(174.3628,659.6893) -- (165.2100,666.3390)
  		(170.8668,627.4110) -- (159.5532,627.4110)
  		(159.5532,627.4110) -- (163.0492,638.1704)
  		(159.5532,648.9298) -- (163.0492,638.1704)
  		(174.3628,659.6893) -- (163.0492,659.6893)
  		(174.3628,659.6893) -- (170.8668,648.9298)
  		(150.4012,655.5795) -- (159.5532,648.9298)
  		(161.7140,620.7613) -- (150.4012,620.7613)
  		(161.7140,620.7613) -- (165.2100,610.0019)
  		(165.2100,610.0019) -- (153.8964,610.0019)
  		(153.8964,610.0019) -- (150.4012,620.7613)
  		(153.8964,610.0019) -- (163.0492,603.3522)
  		(141.2484,627.4110) -- (150.4012,620.7613)
  		(139.0876,620.7613) -- (150.4012,620.7613)
  		(141.2484,614.1115) -- (150.4012,620.7613)
  		(141.2484,614.1115) -- (144.7444,603.3522)
  		(144.7444,603.3522) -- (153.8964,610.0019)
  		(153.8964,610.0019) -- (150.4012,599.2424)
  		(129.9356,614.1115) -- (133.4308,603.3522)
  		(135.5916,610.0019) -- (129.9356,614.1115)
  		(141.2484,614.1115) -- (136.9268,614.1115)
  		(205.8579,620.5829) -- (217.1707,620.5829)
  		(208.0179,613.9332) -- (217.1707,620.5829)
  		(217.1707,620.5829) -- (208.0179,627.2326)
  		(211.5139,603.1738) -- (208.0179,613.9332)
  		(211.5139,637.9921) -- (208.0179,627.2326)
  		(217.1707,620.5829) -- (220.6667,631.3424)
  		(220.6667,631.3424) -- (229.8195,637.9921)
  		(211.5139,637.9921) -- (220.6667,631.3424)
  		(211.5139,637.9921) -- (220.6667,644.6417)
  		(217.1707,620.5829) -- (220.6667,609.8235)
  		(220.6667,609.8235) -- (211.5139,603.1738)
  		(220.6667,609.8235) -- (231.9795,609.8235)
  		(231.9795,609.8235) -- (228.4835,620.5829)
  		(237.6363,627.2326) -- (228.4835,620.5829)
  		(217.1707,620.5829) -- (226.3235,627.2326)
  		(226.3235,627.2326) -- (229.8195,637.9921)
  		(226.3235,627.2326) -- (237.6363,627.2326)
  		(241.1323,637.9921) -- (229.8195,637.9921)
  		(226.3235,648.7515) -- (229.8195,637.9921)
  		(226.3235,648.7515) -- (237.6363,648.7515)
  		(241.1323,637.9921) -- (237.6363,648.7515)
  		(226.3235,648.7515) -- (229.8195,659.5109)
  		(241.1323,659.5109) -- (229.8195,659.5109)
  		(241.1323,659.5109) -- (237.6363,648.7515)
  		(226.3235,648.7515) -- (217.1707,655.4011)
  		(220.6667,666.1606) -- (229.8195,659.5109)
  		(220.6667,666.1606) -- (217.1707,655.4011)
  		(205.8579,655.4011) -- (209.3539,666.1606)
  		(205.8579,655.4011) -- (217.1707,655.4011)
  		(241.1323,659.5109) -- (231.9795,666.1606)
  		(241.1323,659.5109) -- (244.6283,648.7515)
  		(220.6667,609.8235) -- (217.1707,599.0641)
  		(229.8195,603.1738) -- (220.6667,609.8235)
  		(241.1323,659.5109) -- (250.2851,666.1606)
  		(217.1707,655.4011) -- (208.0179,648.7515)
  		(217.1707,655.4011) -- (220.6667,644.6417)
  		(229.8195,637.9921) -- (220.6667,644.6417)
  		(217.1707,620.5829) -- (228.4835,620.5829)
  		(208.0179,613.9332) -- (203.6971,613.9332)
  		(226.3235,592.4144) -- (250.2851,666.1606)
  		(250.2851,666.1606) -- (202.3619,666.1606)
  		(211.5139,637.9921) -- (202.3619,609.8235)
  		(202.3619,609.8235) -- (226.3235,592.4144)
  		(202.3619,666.1606) -- (211.5139,637.9921)
  		(159.5532,592.5927) -- (165.2100,610.0019)
  		(161.7140,620.7613) -- (170.8668,627.4110)
  		(170.8668,627.4110) -- (177.8580,648.9298)
  		(174.3628,659.6893) -- (183.5148,666.3390)
  		(183.5148,666.3390) -- (177.8580,648.9298)
  		(165.2100,610.0019) -- (170.8668,627.4110)
  		(159.5532,592.5927) -- (144.7444,603.3522)
  		(144.7444,603.3522) -- (133.4308,603.3522)
  		(129.9356,614.1115) -- (139.0876,620.7613)
  		(139.0876,620.7613) -- (135.5916,610.0019)
  		(135.5916,610.0019) -- (144.7444,603.3522)
  		(144.7444,638.1704) -- (135.5916,666.3390)
  		(139.0876,655.5795) -- (142.5836,666.3390)
  		(135.5916,666.3390) -- (142.5836,666.3390)
  		(142.5836,666.3390) -- (165.2100,666.3390)
  		(174.3628,672.9887) -- (183.5148,666.3390)
  		(183.5148,666.3390) -- (165.2100,666.3390)
  		(174.3628,672.9887) -- (165.2100,666.3390)
  		(139.0876,620.7613) -- (144.7444,638.1704);
  		
  		\draw[line width=0.8pt]
  		(50.0806,648.8731) -- (46.5846,638.1138)
  		(46.5846,638.1138) -- (43.0886,627.3543)
  		(37.4318,666.2823) -- (26.1190,666.2823)
  		(11.3102,655.5228) -- (14.8054,666.2823)
  		(16.9662,638.1138) -- (13.4702,648.8731)
  		(50.0806,648.8731) -- (46.5846,659.6326)
  		(37.4318,666.2823) -- (46.5846,672.9320)
  		(46.5846,659.6326) -- (55.7374,666.2823)
  		(46.5846,672.9320) -- (55.7374,666.2823)
  		(11.3102,655.5228) -- (14.8054,666.2823)
  		(26.1190,666.2823) -- (14.8054,666.2823)
  		(43.0886,627.3543) -- (33.9358,620.7046)
  		(33.9358,620.7046) -- (37.4318,609.9452)
  		(35.2710,603.2955) -- (31.7758,592.5361)
  		(2.1574,614.0550) -- (11.3102,620.7046)
  		(2.1574,614.0550) -- (5.6534,603.2955)
  		(5.6534,603.2955) -- (16.9662,603.2955)
  		(22.6230,599.1858) -- (31.7758,592.5361)
  		(16.9662,638.1138) -- (13.4702,627.3543)
  		(11.3102,655.5228) -- (13.4702,648.8731)
  		(16.9662,603.2955) -- (22.6230,599.1858)
  		(37.4318,609.9452) -- (35.2710,603.2955)
  		(13.4702,627.3543) -- (11.3102,620.7046)
  		(112.6818,648.9364) -- (109.1858,638.1770)
  		(109.1858,638.1770) -- (105.6898,627.4175)
  		(100.0330,666.3455) -- (88.7202,666.3455)
  		(73.9106,655.5861) -- (77.4066,666.3455)
  		(79.5674,638.1770) -- (76.0714,648.9364)
  		(112.6818,648.9364) -- (109.1858,659.6958)
  		(100.0330,666.3455) -- (109.1858,672.9952)
  		(109.1858,659.6958) -- (118.3378,666.3455)
  		(109.1858,672.9952) -- (118.3378,666.3455)
  		(88.7202,666.3455) -- (77.4066,666.3455)
  		(105.6898,627.4175) -- (96.5370,620.7678)
  		(96.5370,620.7678) -- (100.0330,610.0084)
  		(97.8722,603.3587) -- (94.3762,592.5993)
  		(64.7586,614.1182) -- (73.9106,620.7678)
  		(64.7586,614.1182) -- (68.2546,603.3587)
  		(68.2546,603.3587) -- (79.5674,603.3587)
  		(85.2242,599.2490) -- (94.3762,592.5993)
  		(79.5674,638.1770) -- (76.0714,627.4175)
  		(73.9106,620.7678) -- (76.0714,627.4175)
  		(73.9106,655.5861) -- (76.0714,648.9364)
  		(79.5674,603.3587) -- (85.2242,599.2490)
  		(97.8722,603.3587) -- (100.0330,610.0084)
  		(177.8580,648.9298) -- (174.3628,659.6893)
  		(161.7140,620.7613) -- (165.2100,610.0019)
  		(129.9356,614.1115) -- (133.4308,603.3522)
  		(135.5916,610.0019) -- (129.9356,614.1115)
  		(226.3235,592.4144) -- (250.2851,666.1606)
  		(250.2851,666.1606) -- (202.3619,666.1606)
  		(211.5139,637.9921) -- (202.3619,609.8235)
  		(202.3619,609.8235) -- (226.3235,592.4144)
  		(202.3619,666.1606) -- (211.5139,637.9921)
  		(159.5532,592.5927) -- (165.2100,610.0019)
  		(161.7140,620.7613) -- (170.8668,627.4110)
  		(170.8668,627.4110) -- (177.8580,648.9298)
  		(174.3628,659.6893) -- (183.5148,666.3390)
  		(183.5148,666.3390) -- (177.8580,648.9298)
  		(165.2100,610.0019) -- (170.8668,627.4110)
  		(159.5532,592.5927) -- (144.7444,603.3522)
  		(144.7444,603.3522) -- (133.4308,603.3522)
  		(129.9356,614.1115) -- (139.0876,620.7613)
  		(139.0876,620.7613) -- (135.5916,610.0019)
  		(135.5916,610.0019) -- (144.7444,603.3522)
  		(144.7444,638.1704) -- (135.5916,666.3390)
  		(139.0876,655.5795) -- (142.5836,666.3390)
  		(135.5916,666.3390) -- (142.5836,666.3390)
  		(142.5836,666.3390) -- (165.2100,666.3390)
  		(174.3628,672.9887) -- (183.5148,666.3390)
  		(183.5148,666.3390) -- (165.2100,666.3390)
  		(174.3628,672.9887) -- (165.2100,666.3390)
  		(139.0876,620.7613) -- (144.7444,638.1704);
  		
  		\fill[line width=0.1pt] (22,606) circle[radius=0.05cm];
  		
  		\draw[red, <-] (172,616) to [bend left] (138,600);
  		\draw[red, ->] (184,672) to [bend left] (185,655);
  		\draw[red, ->] (134,620) to [bend left] (133,661);
  		
  		\small
  		\coordinate[label=right:Dog] (a) at (19,680);
  		\coordinate[label=right:$(b)$] (a) at (88,680);
  		\coordinate[label=right:$(c)$] (a) at (153,680);
  		\coordinate[label=right:$T_9$] (a) at (220,680);
  		\normalsize
  	\end{tikzpicture}	
  	\caption{Building $T_9$ from the Dog.}
  \end{figure}
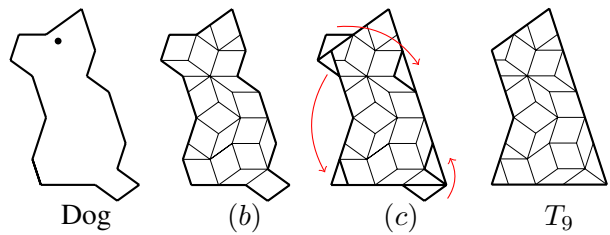
  
  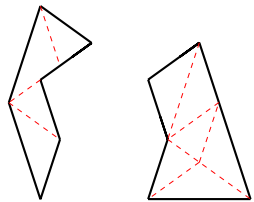
\begin{figure}[!ht]
  		\centering
  		\begin{tikzpicture}
  		[y=-1.0pt, x=1.0pt, scale=0.8]
  		\draw[line width=0.8pt]
  		(31.0652,605.3026) -- (40.2180,598.6529)
  		(31.0652,605.3026) -- (40.2180,598.6529)
  		(16.2564,616.0619) -- (40.2180,598.6529)
  		(75.8175,644.1336) -- (72.3215,633.3736)
  		(81.4743,605.2054) -- (90.6271,598.5557)
  		(81.4743,605.2054) -- (90.6271,598.5557)
  		(90.6271,598.5557) -- (114.5887,672.3016)
  		(114.5887,672.3016) -- (66.6655,672.3016)
  		(75.8175,644.1336) -- (66.6655,615.9648)
  		(66.6655,615.9648) -- (90.6271,598.5557)
  		(16.2564,581.2438) -- (40.2180,598.6529)
  		(1.4476,626.8211) -- (16.2564,672.3987)
  		(16.2564,616.0619) -- (25.4092,644.2307)
  		(1.4476,626.8211) -- (16.2564,581.2438)
  		(16.2564,672.3987) -- (25.4092,644.2307)
  		(75.8175,644.1336) -- (66.6655,672.3016);
  		
  		\draw[red, dash pattern=on 2pt off 2pt, line width=0.2pt]
  		(16.2564,616.0619) -- (1.4476,626.8211)
  		(1.4476,626.8211) -- (25.4092,644.2307)
  		(16.2564,581.2438) -- (25.4092,609.4123)
  		(90.6271,598.5557) -- (75.8175,644.1336)
  		(75.8175,644.1336) -- (114.5887,672.3016)
  		(90.6271,654.8928) -- (66.6655,672.3016)
  		(90.6271,654.8928) -- (99.7791,626.7240)
  		(99.7791,626.7240) -- (75.8175,644.1336);
  		\end{tikzpicture}
  	\caption{$T_8$ and $T_9$ with Robinson triangles.}
  \end{figure}
  
  As well as the Penrose tiles, $T_8$ and $T_9$ can be decomposed into Robinson triangles. $T_8$ can be decomposed into two golden gnomons and two golden triangles, $T_9$ into four golden gnomons and two golden triangles (Fig. 7).
  
  The decorated snake and dog tiles (Fig. 5$b$ and 6$b$) lead to a P3 tiling, whereas a tiling of the decorated $T_8$ and  $T_9$ tiles will always contain small errors due to incomplete edges. These incorrect areas are related to the decoration in the top corner of $T_8$ and the lower left corner of $T_9$. Figure 8 shows an aperiodic set of two decorated prototiles (based on $T_8$ and $T_9$), which allow an exact P3 tiling again (Fig. 12). We leave it to the reader to locate the three areas with the rectified edge-errors.
  
  \begin{figure}[!ht]
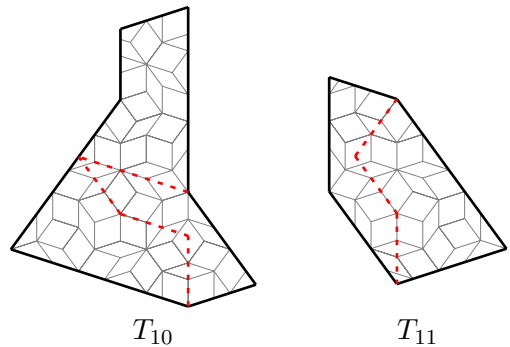

  	\centering

  	\caption{Aperiodic set for a P3 tiling.}
  \end{figure}
  
  Without substitution rules the undecorated $T_8$, $T_9$, and $T_{11}$ tiles also allow periodic tilings (Fig. 9). A tessellation of the plane only with $T_8$ or $T_{10}$ tiles is not possible.
  
  Figures 11 and 12 show larger portions with the snake and dog tiles and the undecorated $T_{10},T_{11}$ set. Note the special feature in the last one, where the $T_{11}$ tiles never touch each other.\\
  
  \begin{figure}[!ht]
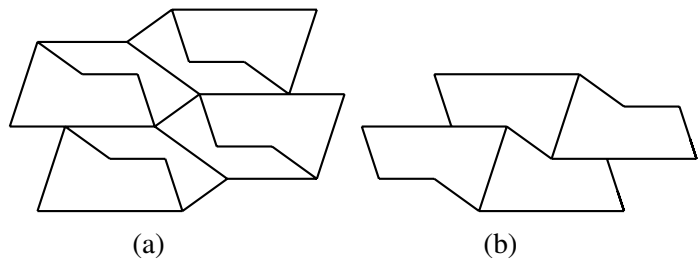

  	\centering
  	\begin{minipage}[t]{0.30\linewidth}
  		\centering
  		%
  		\subcaption{}
  	\end{minipage}
  	\quad\quad
  	\begin{minipage}[t]{0.30\linewidth}
  		\centering
  		%
  		\subcaption{}
  	\end{minipage}
  	\caption{Matching rules for periodic tilings with $T_8$, $T_9$ and $T_{11}$.}
  \end{figure}
  
  \newpage
  
  The substitution rules for all aperiodic sets in this article are the same as for the Penrose kite and dart, due to their same mld-class. Figure 10 shows these rules using $T_8$ and $T_9$ as an example.
  
  \begin{figure}[!ht]
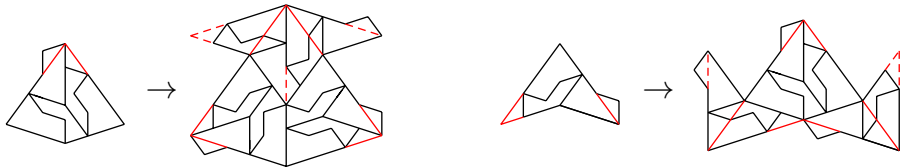

  	\centering
  	%
  	\caption{Substitution rules for $T_8$ and $T_9$.}
  \end{figure}
  
  \begin{figure}[!ht]
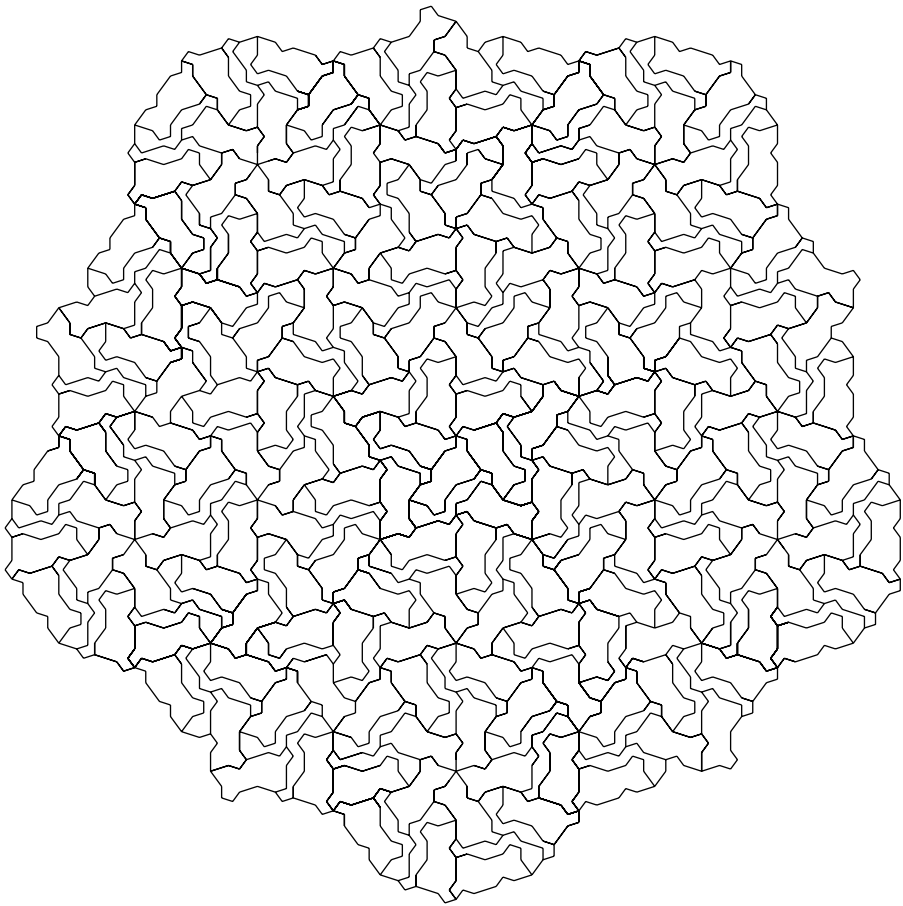

	\centering
	%
	\caption{Snake and Dog tiling.}
  \end{figure}
  
  \newpage
  
  \begin{figure}[!ht]
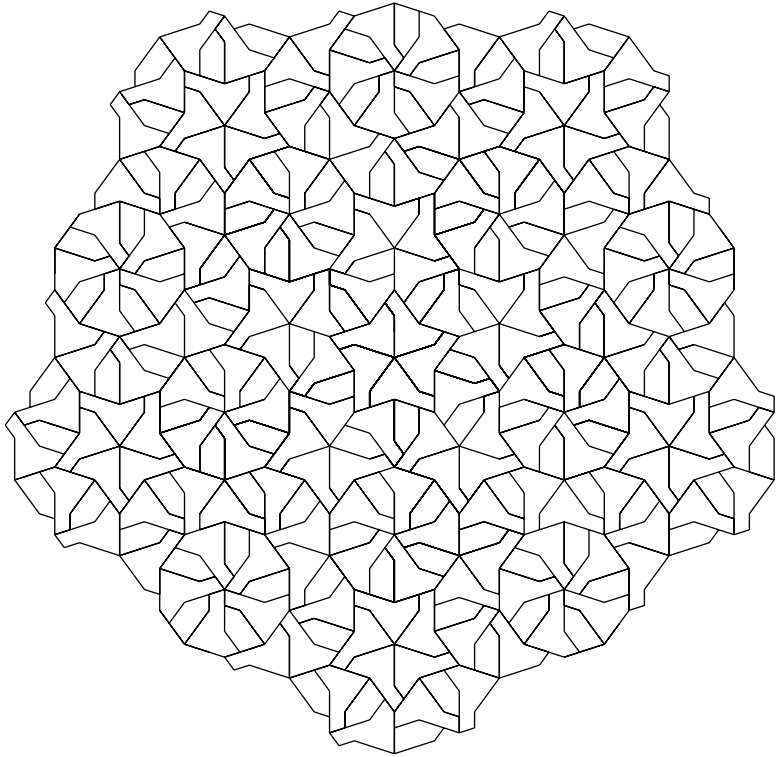

  	\centering
  	%
  \caption{Tiling with the undecorated $T_{10},T_{11}$ set.}
  \end{figure}

  \section{Aperiodic Monotiles} 
  
  One of the most exciting open problems in plane geometry is the existence of an aperiodic monotile. It asks about a single connected prototile that by itself forms a strongly aperiodic set. Such a tile can tessellate the Euclidean plane only in a nonperiodic way without matching rules. The smallest aperiodic sets known to date consist of two prototiles, like the Penrose tiles. Note that the P2 and P3 sets must be modified to get strongly aperiodic sets without matching rules. Examples can be seen in \cite{Gruenbaum} on page 539 (kite, dart) and page 544 (rhombs).\footnote{Further examples: https://en.wikipedia.org/wiki/Penrose\_tiling\#Rhombus\_tiling\_(P3)}
  
  The currently best approximations to an aperiodic monotile were given by Petra Gummelt in 1996, and  Joshua Socolar and Joan Taylor in 2010. Gummelt constructed a decorated decagonal tile and showed that when two kinds of \textit{overlaps} between pairs of tiles are allowed, these tiles can cover the plane only nonperiodically \cite{Gummelt}. Socolar and Taylor presented an undecorated, but \textit{not connected} aperiodic monotile that is based on a regular hexagon \cite{Socolar}.
  
  Figure 13 shows an approximation to an aperiodic monotile from the Penrose tiles, a shape that represents both kite and dart as well as possible. Please note that it is not possible to create a \textit{truly} aperiodic monotile with a 5-fold rotational symmetry from the Penrose prototiles because they do not have the same acreage. But other shapes as given by $T_{12}$ are possible. It is also possible to modify edges (e.g. by identically shaped bumps and notches) to enforce the tiling rules.
  
  $T_{12}$, an irregular concave dodecagon, is based on a connected pair of $T_{10},T_{11}$ tiles (Fig. 13a), and can be subdivided into a rhombus $(a,b,c,p)$, two congruent kites $(c,d,r,p)$ and $(g,h,i,r)$, two non-congruent trapezoids $(d,e,f,r)$ and $(i,j,k,r)$, which we will denote by $T_{13}$ and $T_{14}$, and an irregular quadrilateral $(l,n,q,r)$ (Fig. 13b).
  
  \begin{figure}[!ht]
  	\centering
  	\begin{minipage}[t]{0.40\linewidth}
  		\centering
  		\begin{tikzpicture}
  			[y=-1.0pt, x=1.0pt, scale=1.0]
  			\draw[black!50, line width=0.01pt]
  			(206.6460,627.7048) -- (206.6460,616.3917)
  			(206.6460,627.7048) -- (213.2956,636.8573)
  			(217.4052,631.2007) -- (224.0556,640.3532)
  			(224.0556,640.3532) -- (213.2956,636.8573)
  			(224.0556,610.7350) -- (234.8148,607.2391)
  			(224.0556,603.7432) -- (234.8148,607.2391)
  			(234.8148,607.2391) -- (228.1652,616.3917)
  			(206.6460,598.0866) -- (206.6460,586.7735)
  			(206.6460,586.7735) -- (213.2956,595.9260)
  			(213.2956,595.9260) -- (224.0556,592.4300)
  			(224.0556,592.4300) -- (224.0556,603.7432)
  			(206.6460,586.7735) -- (217.4052,583.2775)
  			(206.6460,586.7735) -- (206.6460,575.4603)
  			(206.6460,586.7735) -- (206.6460,575.4603)
  			(234.8148,625.5442) -- (228.1652,616.3917)
  			(234.8148,607.2391) -- (241.4644,616.3917)
  			(241.4644,616.3917) -- (252.2236,619.8876)
  			(234.8148,625.5442) -- (241.4644,616.3917)
  			(234.8148,625.5442) -- (245.5740,629.0401)
  			(234.8148,607.2391) -- (234.8148,595.9260)
  			(234.8148,595.9260) -- (224.0556,592.4300)
  			(234.8148,595.9260) -- (245.5740,592.4300)
  			(245.5740,592.4300) -- (245.5740,603.7432)
  			(256.3332,607.2391) -- (245.5740,603.7432)
  			(234.8148,607.2391) -- (245.5740,610.7350)
  			(245.5740,610.7350) -- (252.2236,619.8876)
  			(245.5740,610.7350) -- (256.3332,607.2391)
  			(262.9836,616.3917) -- (252.2236,619.8876)
  			(262.9836,616.3917) -- (256.3332,607.2391)
  			(252.2236,631.2007) -- (252.2236,619.8876)
  			(252.2236,631.2007) -- (262.9836,627.7048)
  			(262.9836,616.3917) -- (252.2236,619.8876)
  			(262.9836,616.3917) -- (262.9836,627.7048)
  			(252.2236,631.2007) -- (258.8732,640.3532)
  			(252.2236,631.2007) -- (262.9836,627.7048)
  			(269.6332,636.8573) -- (258.8732,640.3532)
  			(269.6332,636.8573) -- (262.9836,627.7048)
  			(252.2236,631.2007) -- (258.8732,640.3532)
  			(252.2236,631.2007) -- (245.5740,640.3532)
  			(252.2236,649.5062) -- (258.8732,640.3532)
  			(252.2236,649.5062) -- (245.5740,640.3532)
  			(234.8148,643.8492) -- (241.4644,653.0014)
  			(234.8148,643.8492) -- (245.5740,640.3532)
  			(252.2236,649.5062) -- (241.4644,653.0014)
  			(252.2236,649.5062) -- (245.5740,640.3532)
  			(234.8148,643.8492) -- (241.4644,653.0014)
  			(234.8148,643.8492) -- (224.0556,640.3532)
  			(224.0556,640.3532) -- (213.2956,636.8573)
  			(241.4644,653.0014) -- (230.7052,649.5062)
  			(230.7052,649.5062) -- (224.0556,640.3532)
  			(234.8148,643.8492) -- (241.4644,653.0014)
  			(234.8148,643.8492) -- (224.0556,640.3532)
  			(206.6460,627.7048) -- (206.6460,616.3917)
  			(206.6460,598.0866) -- (206.6460,586.7735)
  			(206.6460,586.7735) -- (206.6460,575.4603)
  			(206.6460,586.7735) -- (206.6460,575.4603)
  			(206.6460,646.0102) -- (213.2956,636.8573)
  			(224.0556,640.3532) -- (213.2956,636.8573)
  			(224.0556,640.3532) -- (217.4052,649.5062)
  			(206.6460,646.0102) -- (213.2956,636.8573)
  			(206.6460,646.0102) -- (217.4052,649.5062)
  			(224.0556,640.3532) -- (230.7052,649.5062)
  			(224.0556,640.3532) -- (217.4052,649.5062)
  			(224.0556,658.6582) -- (230.7052,649.5062)
  			(224.0556,658.6582) -- (217.4052,649.5062)
  			(241.4644,653.0014) -- (230.7052,649.5062)
  			(224.0556,658.6582) -- (230.7052,649.5062)
  			(206.6460,646.0102) -- (213.2956,636.8573)
  			(206.6460,653.0014) -- (217.4052,649.5062)
  			(206.6460,653.0014) -- (206.6460,664.3150)
  			(206.6460,664.3150) -- (217.4052,660.8190)
  			(217.4052,660.8190) -- (217.4052,649.5062)
  			(224.0556,658.6582) -- (217.4052,660.8190)
  			(269.6332,636.8573) -- (262.9836,646.0102)
  			(252.2236,649.5062) -- (262.9836,646.0102)
  			(269.6332,636.8573) -- (269.6332,625.5442)
  			(269.6332,625.5442) -- (262.9836,616.3917)
  			(228.1652,586.7735) -- (217.4052,583.2775)
  			(234.8148,595.9260) -- (228.1652,586.7735)
  			(234.8148,595.9260) -- (228.1652,586.7735)
  			(228.1652,586.7735) -- (234.8148,577.6210)
  			(234.8148,577.6210) -- (241.4644,586.7735)
  			(241.4644,586.7735) -- (234.8148,595.9260)
  			(217.4052,583.2775) -- (228.1652,586.7735)
  			(217.4052,583.2775) -- (228.1652,586.7735)
  			(228.1652,586.7735) -- (234.8148,577.6210)
  			(234.8148,577.6210) -- (224.0556,574.1250)
  			(224.0556,574.1250) -- (217.4052,583.2775)
  			(234.8148,577.6210) -- (224.0556,574.1250)
  			(224.0556,658.6582) -- (241.4644,653.0014)
  			(256.3332,607.2391) -- (245.5740,592.4300)
  			(245.5740,592.4300) -- (241.4644,586.7735)
  			(269.6332,636.8573) -- (280.3924,640.3532)
  			(280.3924,640.3532) -- (269.6332,625.5442)
  			(262.9836,646.0102) -- (280.3924,640.3532)
  			(206.6460,627.7048) -- (206.6460,616.3917)
  			(206.6460,627.7048) -- (199.9964,636.8573)
  			(195.8868,631.2007) -- (189.2372,640.3532)
  			(195.8868,619.8876) -- (189.2372,629.0401)
  			(189.2372,629.0401) -- (189.2372,640.3532)
  			(189.2372,640.3532) -- (199.9964,636.8573)
  			(195.8868,619.8876) -- (185.1276,616.3917)
  			(189.2372,610.7350) -- (178.4772,607.2391)
  			(189.2372,603.7432) -- (178.4772,607.2391)
  			(178.4772,607.2391) -- (185.1276,616.3917)
  			(206.6460,598.0866) -- (206.6460,586.7735)
  			(206.6460,586.7735) -- (199.9964,595.9260)
  			(199.9964,595.9260) -- (189.2372,592.4300)
  			(189.2372,592.4300) -- (189.2372,603.7432)
  			(206.6460,586.7735) -- (195.8868,583.2775)
  			(206.6460,586.7735) -- (206.6460,575.4603)
  			(206.6460,586.7735) -- (206.6460,575.4603)
  			(178.4772,625.5442) -- (189.2372,629.0401)
  			(178.4772,625.5442) -- (185.1276,616.3917)
  			(178.4772,607.2391) -- (171.8276,616.3917)
  			(171.8276,616.3917) -- (161.0684,619.8876)
  			(178.4772,607.2391) -- (185.1276,616.3917)
  			(178.4772,625.5442) -- (171.8276,616.3917)
  			(178.4772,625.5442) -- (185.1276,616.3917)
  			(161.0684,619.8876) -- (171.8276,616.3917)
  			(171.8276,616.3917) -- (178.4772,607.2391)
  			(161.0684,619.8876) -- (167.7180,629.0401)
  			(178.4772,625.5442) -- (171.8276,616.3917)
  			(178.4772,625.5442) -- (167.7180,629.0401)
  			(167.7180,640.3532) -- (178.4772,636.8573)
  			(178.4772,636.8573) -- (189.2372,640.3532)
  			(167.7180,640.3532) -- (167.7180,629.0401)
  			(178.4772,625.5442) -- (178.4772,636.8573)
  			(178.4772,625.5442) -- (167.7180,629.0401)
  			(178.4772,607.2391) -- (178.4772,595.9260)
  			(178.4772,595.9260) -- (189.2372,592.4300)
  			(178.4772,595.9260) -- (167.7180,592.4300)
  			(167.7180,592.4300) -- (167.7180,603.7432)
  			(156.9588,607.2391) -- (167.7180,603.7432)
  			(178.4772,607.2391) -- (167.7180,610.7350)
  			(167.7180,610.7350) -- (161.0684,619.8876)
  			(167.7180,610.7350) -- (156.9588,607.2391)
  			(150.3092,616.3917) -- (161.0684,619.8876)
  			(150.3092,616.3917) -- (156.9588,607.2391)
  			(161.0684,631.2007) -- (161.0684,619.8876)
  			(161.0684,631.2007) -- (150.3092,627.7048)
  			(150.3092,616.3917) -- (161.0684,619.8876)
  			(150.3092,616.3917) -- (150.3092,627.7048)
  			(161.0684,631.2007) -- (154.4188,640.3532)
  			(161.0684,631.2007) -- (150.3092,627.7048)
  			(143.6596,636.8573) -- (154.4188,640.3532)
  			(143.6596,636.8573) -- (150.3092,627.7048)
  			(161.0684,631.2007) -- (154.4188,640.3532)
  			(161.0684,631.2007) -- (167.7180,640.3532)
  			(161.0684,649.5062) -- (154.4188,640.3532)
  			(161.0684,649.5062) -- (167.7180,640.3532)
  			(178.4772,643.8492) -- (171.8276,653.0014)
  			(178.4772,643.8492) -- (167.7180,640.3532)
  			(161.0684,649.5062) -- (171.8276,653.0014)
  			(161.0684,649.5062) -- (167.7180,640.3532)
  			(178.4772,643.8492) -- (171.8276,653.0014)
  			(178.4772,643.8492) -- (189.2372,640.3532)
  			(189.2372,640.3532) -- (199.9964,636.8573)
  			(171.8276,653.0014) -- (182.5876,649.5062)
  			(182.5876,649.5062) -- (189.2372,640.3532)
  			(178.4772,643.8492) -- (171.8276,653.0014)
  			(178.4772,643.8492) -- (189.2372,640.3532)
  			(206.6460,598.0866) -- (206.6460,586.7735)
  			(206.6460,586.7735) -- (206.6460,575.4603)
  			(206.6460,646.0102) -- (199.9964,636.8573)
  			(189.2372,640.3532) -- (199.9964,636.8573)
  			(189.2372,640.3532) -- (195.8868,649.5062)
  			(206.6460,646.0102) -- (199.9964,636.8573)
  			(206.6460,646.0102) -- (195.8868,649.5062)
  			(189.2372,640.3532) -- (182.5876,649.5062)
  			(189.2372,640.3532) -- (195.8868,649.5062)
  			(189.2372,658.6582) -- (182.5876,649.5062)
  			(189.2372,658.6582) -- (195.8868,649.5062)
  			(171.8276,653.0014) -- (182.5876,649.5062)
  			(189.2372,658.6582) -- (182.5876,649.5062)
  			(206.6460,646.0102) -- (199.9964,636.8573)
  			(206.6460,653.0014) -- (195.8868,649.5062)
  			(206.6460,653.0014) -- (206.6460,664.3150)
  			(206.6460,664.3150) -- (195.8868,660.8190)
  			(195.8868,660.8190) -- (195.8868,649.5062)
  			(189.2372,658.6582) -- (195.8868,660.8190)
  			(143.6596,636.8573) -- (150.3092,646.0102)
  			(161.0684,649.5062) -- (150.3092,646.0102)
  			(143.6596,636.8573) -- (143.6596,625.5442)
  			(143.6596,625.5442) -- (150.3092,616.3917)
  			(185.1276,586.7735) -- (195.8868,583.2775)
  			(178.4772,595.9260) -- (185.1276,586.7735)
  			(178.4772,595.9260) -- (185.1276,586.7735)
  			(185.1276,586.7735) -- (178.4772,577.6210)
  			(178.4772,577.6210) -- (171.8276,586.7735)
  			(171.8276,586.7735) -- (178.4772,595.9260)
  			(195.8868,583.2775) -- (185.1276,586.7735)
  			(195.8868,583.2775) -- (185.1276,586.7735)
  			(185.1276,586.7735) -- (178.4772,577.6210)
  			(178.4772,577.6210) -- (189.2372,574.1250)
  			(189.2372,574.1250) -- (195.8868,583.2775)
  			(178.4772,577.6210) -- (189.2372,574.1250)
  			(189.2372,574.1250) -- (189.2372,562.8119)
  			(195.8868,571.9644) -- (189.2372,562.8119)
  			(189.2372,658.6582) -- (171.8276,653.0014)
  			(156.9588,607.2391) -- (167.7180,592.4300)
  			(167.7180,592.4300) -- (171.8276,586.7735)
  			(143.6596,636.8573) -- (132.8996,640.3532)
  			(132.8996,640.3532) -- (143.6596,625.5442)
  			(150.3092,646.0102) -- (132.8996,640.3532)
  			(206.6460,598.0866) -- (213.2956,607.2391)
  			(206.6460,598.0866) -- (199.9964,607.2391)
  			(199.9964,607.2391) -- (189.2372,610.7350)
  			(189.2372,610.7350) -- (195.8868,619.8876)
  			(199.9964,607.2391) -- (206.6460,616.3917)
  			(195.8868,631.2007) -- (206.6460,627.7048)
  			(195.8868,619.8876) -- (206.6460,616.3917)
  			(195.8868,619.8876) -- (195.8868,631.2007)
  			(224.0556,610.7350) -- (217.4052,619.8876)
  			(213.2956,607.2391) -- (224.0556,610.7350)
  			(213.2956,607.2391) -- (206.6460,616.3917)
  			(217.4052,619.8876) -- (206.6460,616.3917)
  			(206.6460,627.7048) -- (206.6460,616.3917)
  			(217.4052,631.2007) -- (206.6460,627.7048)
  			(217.4052,619.8876) -- (217.4052,631.2007)
  			(217.4052,619.8876) -- (228.1652,616.3917)
  			(217.4052,619.8876) -- (224.0556,629.0401)
  			(241.4644,616.3917) -- (234.8148,607.2391)
  			(234.8148,607.2391) -- (228.1652,616.3917)
  			(245.5740,640.3532) -- (234.8148,636.8573)
  			(245.5740,640.3532) -- (245.5740,629.0401)
  			(234.8148,625.5442) -- (224.0556,629.0401)
  			(224.0556,629.0401) -- (224.0556,640.3532)
  			(234.8148,636.8573) -- (224.0556,640.3532)
  			(234.8148,625.5442) -- (234.8148,636.8573)
  			(234.8148,625.5442) -- (241.4644,616.3917)
  			(234.8148,625.5442) -- (245.5740,629.0401)
  			(252.2236,619.8876) -- (245.5740,629.0401)
  			(252.2236,619.8876) -- (241.4644,616.3917)
  			(224.0556,592.4300) -- (217.4052,583.2775)
  			(206.6460,586.7735) -- (213.2956,595.9260)
  			(213.2956,595.9260) -- (224.0556,592.4300)
  			(213.2956,607.2391) -- (213.2956,595.9260)
  			(213.2956,607.2391) -- (224.0556,603.7432)
  			(234.8148,595.9260) -- (224.0556,592.4300)
  			(224.0556,592.4300) -- (224.0556,603.7432)
  			(234.8148,607.2391) -- (224.0556,603.7432)
  			(206.6460,575.4603) -- (217.4052,571.9644)
  			(206.6460,586.7735) -- (217.4052,583.2775)
  			(217.4052,571.9644) -- (217.4052,583.2775)
  			(245.5740,592.4300) -- (245.5740,603.7432)
  			(234.8148,595.9260) -- (245.5740,592.4300)
  			(234.8148,607.2391) -- (234.8148,595.9260)
  			(234.8148,607.2391) -- (245.5740,603.7432)
  			(206.6460,586.7735) -- (206.6460,575.4603)
  			(206.6460,575.4603) -- (195.8868,571.9644)
  			(195.8868,571.9644) -- (195.8868,583.2775)
  			(189.2372,592.4300) -- (195.8868,583.2775)
  			(206.6460,586.7735) -- (199.9964,595.9260)
  			(206.6460,586.7735) -- (195.8868,583.2775)
  			(199.9964,595.9260) -- (189.2372,592.4300)
  			(199.9964,607.2391) -- (189.2372,603.7432)
  			(199.9964,607.2391) -- (199.9964,595.9260)
  			(189.2372,592.4300) -- (189.2372,603.7432)
  			(178.4772,595.9260) -- (189.2372,592.4300)
  			(178.4772,607.2391) -- (189.2372,603.7432)
  			(178.4772,607.2391) -- (178.4772,595.9260)
  			(178.4772,595.9260) -- (167.7180,592.4300)
  			(167.7180,592.4300) -- (167.7180,603.7432)
  			(178.4772,607.2391) -- (167.7180,603.7432)
  			(206.6460,575.4603) -- (199.9964,566.3078)
  			(199.9964,566.3078) -- (189.2372,562.8119)
  			(199.9964,566.3078) -- (206.6460,557.1553)
  			(206.6460,557.1553) -- (195.8868,553.6594)
  			(195.8868,553.6594) -- (189.2372,562.8119)
  			(195.8868,553.6594) -- (206.6460,550.1634)
  			(206.6460,550.1634) -- (206.6460,538.8503)
  			(178.4772,577.6210) -- (178.4772,566.3078)
  			(178.4772,577.6210) -- (178.4772,566.3078)
  			(178.4772,566.3078) -- (189.2372,562.8119)
  			(178.4772,559.3159) -- (189.2372,562.8119)
  			(182.5876,553.6594) -- (189.2372,562.8119)
  			(182.5876,553.6594) -- (189.2372,544.5069)
  			(189.2372,544.5069) -- (195.8868,553.6594)
  			(195.8868,553.6594) -- (189.2372,562.8119)
  			(195.8868,542.3462) -- (195.8868,553.6594)
  			(182.5876,553.6594) -- (178.4772,552.3240)
  			(178.4780,655.1622) -- (189.2372,640.3532);
  			
  			\draw[line width=1.0pt]
  			(178.4772,577.6210) -- (178.4772,548.0028)
  			(178.4772,548.0028) -- (206.6460,538.8503)
  			(206.6460,664.3150) -- (178.4780,655.1622)
  			(178.4780,655.1622) -- (189.2372,640.3532)
  			(206.6460,664.3150) -- (252.2236,649.5062)
  			(189.2372,640.3532) -- (147.7676,619.8894)
  			(147.7676,619.8894) -- (178.4772,577.6210)
  			(206.6460,586.7735) -- (237.6068,581.4645)
  			(252.2236,601.5826) -- (237.6068,581.4645)
  			(230.7052,608.5745) -- (252.2236,649.5062)
  			(206.6460,538.8503) -- (206.6460,586.7735)
  			(230.7052,608.5745) -- (252.2236,601.5826);
  			
  			\draw[red, dash pattern=on 2pt off 3pt, line width=1.0]
  			(178.4772,577.6210) -- (206.6460,538.8503)
  			(234.8148,577.6210) -- (206.6460,568.4685)
  			(178.4780,655.1622) -- (132.8996,640.3532)
  			(132.8996,640.3532) -- (147.7676,619.8894)
  			(252.2236,649.5062) -- (280.3924,640.3532)
  			(280.3924,640.3532) -- (252.2236,601.5826)
  			(237.6068,581.4645) -- (234.8148,577.6210);
  		\end{tikzpicture}
  		\subcaption{}
  	\end{minipage}
  	\quad\quad
  	\begin{minipage}[t]{0.40\linewidth}
  		\centering
  		\begin{tikzpicture}
  			[y=-1.0pt, x=1.0pt, xscale=1.0]
  			\draw[line width=1.0pt]
  			(40.4026,546.2207) -- (68.5714,537.0679)
  			(40.4026,546.2207) -- (40.4026,575.8383)
  			(68.5714,537.0679) -- (68.5714,566.6863)
  			(68.5714,566.6863) -- (68.5714,584.9911)
  			(68.5714,662.5319) -- (40.4026,653.3799)
  			(40.4026,653.3799) -- (51.1618,638.5711)
  			(40.4026,575.8383) -- (22.9938,599.7999)
  			(22.9938,599.7999) -- (9.6946,618.1047)
  			(9.6946,618.1047) -- (51.1618,638.5711)
  			(68.5714,662.5319) -- (96.7394,653.3799)
  			(96.7394,653.3799) -- (114.1482,647.7239)
  			(114.1482,647.7239) -- (92.6298,606.7919)
  		    (92.6298,606.7919) -- (114.1482,599.7999)
  			(99.4875,579.7020) -- (68.5714,584.9911)
  			(99.4875,579.7020) -- (114.1490,599.8007);
  			
  			\draw[red, line width=0.1]
  			(40.4026,575.8383) -- (68.5714,566.6863)
  			(68.5714,614.6095) -- (68.5714,584.9911)
  			(68.5714,614.6095) -- (51.1618,638.5711)
  			(68.5714,614.6095) -- (96.7394,653.3799)
  			(68.5714,614.6095) -- (22.9938,599.7999)
  			(68.5714,614.6095) -- (92.6298,606.7919);
  			
  			\draw[red, line width=0.1]
  			(96.7394,575.8383) -- (68.5714,566.6863)
  			(96.7394,575.8383) -- (99.4875,579.7020)
  		    (99.4875,579.7020) -- (114.1490,577.1743)
  			(114.1490,599.8007) -- (114.1490,577.1743);
  			
  			\small
  			\coordinate[label=right:$a$] (a) at (67,534);
  			\coordinate[label=right:$b$] (b) at (28,545);
  			\coordinate[label=right:$c$] (c) at (28,575);
  			\coordinate[label=right:$d$] (d) at (10,597);
  			\coordinate[label=right:$e$] (e) at (-2,618);
  			\coordinate[label=right:$f$] (f) at (47,641);
  			\coordinate[label=right:$g$] (g) at (29,655);
  			\coordinate[label=right:$h$] (h) at (63,655);
  			\coordinate[label=right:$i$] (i) at (92,647);
  			\coordinate[label=right:$j$] (j) at (113,648);
  			\coordinate[label=right:$k$] (k) at (81,603);
  			\coordinate[label=right:$l$] (l) at (114,600);
  			\coordinate[label=right:$m$] (m) at (114,576);
  			\coordinate[label=right:$n$] (n) at (90,587);
  			\coordinate[label=right:$o$] (o) at (94,572);
  			\coordinate[label=right:$p$] (p) at (67,562);
  			\coordinate[label=right:$q$] (q) at (57,585);
  			\coordinate[label=right:$r$] (r) at (57,613);
  			
  			\footnotesize
  			\coordinate[label=right:$T_{13}$] (c) at (26,618);
  			\coordinate[label=right:$T_{14}$] (e) at (80,628);
  		\end{tikzpicture}
  		\subcaption{}
  	\end{minipage}
  	\caption{The approximated aperiodic monotile $T_{12}$.}
  \end{figure}
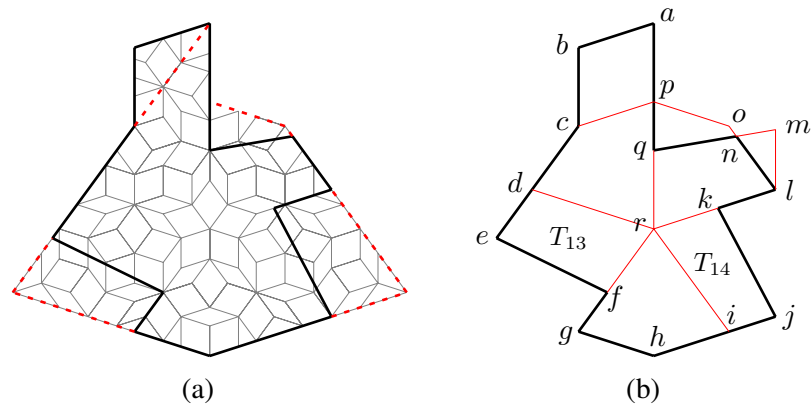
  
  Several properties and common features of the Penrose tilings involve the golden ratio $\varphi=(1+\sqrt{5})/2\approx1.618$ \cite{Gruenbaum}. The edge lengths in Figure 13b, in relation to the unit length edges of the decorated rhombs, are $\vert fg \vert=\vert ij \vert=\vert pq \vert=\varphi$, $\vert ab \vert=\vert bc \vert=\vert cd \vert=\vert gh \vert=\vert hi \vert=\vert lo \vert=\vert op \vert=\vert pa \vert=\vert pc \vert=\vert qr \vert=\vert fr \vert=\varphi+1$, $\vert dr \vert=\vert gr \vert=\vert hj \vert=\vert ir \vert=\vert lr \vert=\vert pr \vert=2\varphi+1$, and $\vert de \vert=\vert kl \vert=\vert lm \vert=2$.
  
  \newpage
  
  A $T_{12}$ tiling contains always five types of gaps. An irregular triangle $(l,m,n)$, an irregular quadrilateral $(n,o,p,q)$, and three types of \textit{rotors} (concave pentadecagons) that can be subdivided into $T_{13}$ and $T_{14}$ tiles (Fig. 13b and 15). Note that a connected pair of $T_{13}$ and $T_{14}$ could build a rhombus with edge lengths of $2\varphi+1$. The decorated $T_{12}$ tiles (Fig. 13a) also allow a P3 tiling with the mentioned gaps. The substitution rules for $T_{12}$ are given in Figure 16. 
  
  \begin{figure}[!ht]
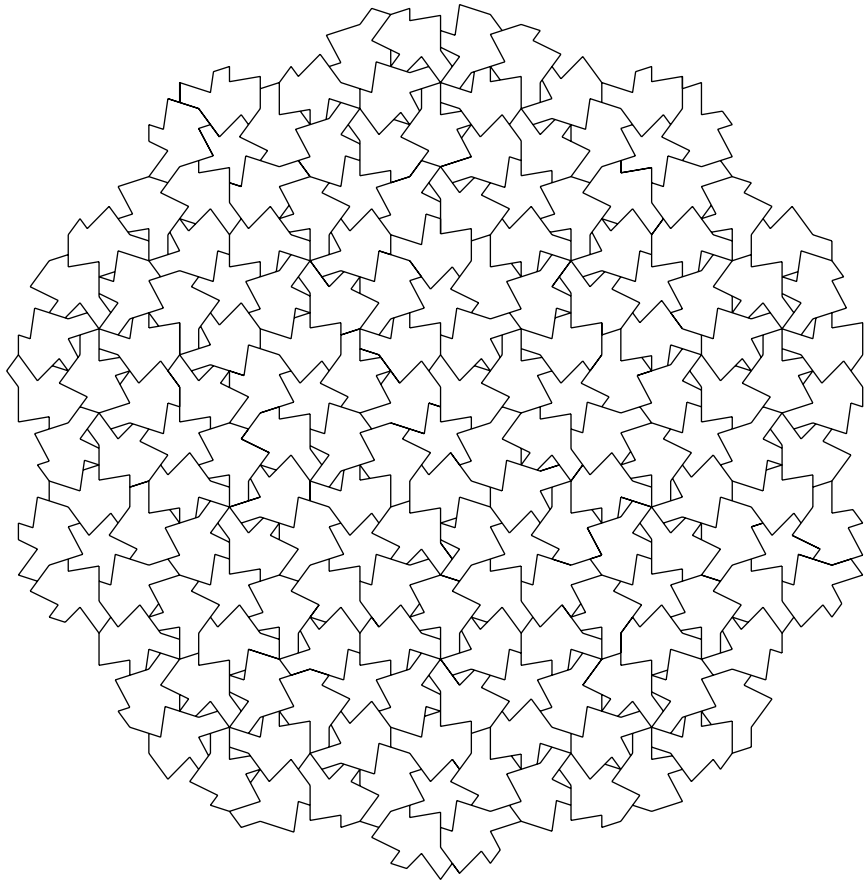

  	\centering

  	\caption{A $T_{12}$ tiling with gaps.}
  \end{figure}
  
  The possibility of an exact P3 tiling by the decorated tile sets mentioned in this article can be proved with the seven possible vertex figures in a P2 tiling. Details are given in \cite{Gruenbaum} on page 561 or at Wikipedia\footnote{https://en.wikipedia.org/wiki/Penrose\_tiling\#Kite\_and\_dart\_tiling\_(P2)}.
  
  \begin{figure}[!ht]
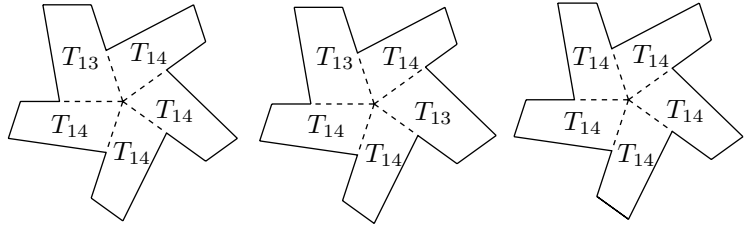

  	\centering
  	%
  	\caption{The three types of \textit{rotor} gaps in a $T_{12}$ tiling.}
  \end{figure}
  \quad
  \begin{figure}[!ht]
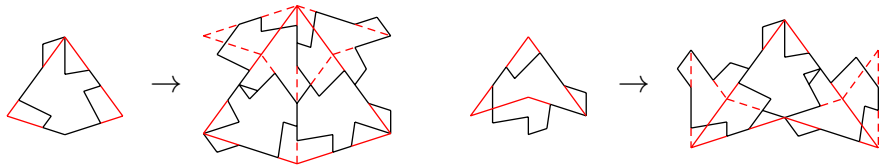

  	\centering
  	%
  	\caption{Substitution rules for $T_{12}$.}
  \end{figure}

\end{document}